# Incircular nets and confocal conics

Arseniy V. Akopyan* Alexander I. Bobenko†

**Abstract**

We consider congruences of straight lines in a plane with the combinatorics of the square grid, with all elementary quadrilaterals possessing an incircle. It is shown that all the vertices of such nets (we call them incircular or IC-nets) lie on confocal conics.

Our main new results are on checkerboard IC-nets in the plane. These are congruences of straight lines in the plane with the combinatorics of the square grid, combinatorially colored as a checkerboard, such that all black coordinate quadrilaterals possess inscribed circles. We show how this larger class of IC-nets appears quite naturally in Laguerre geometry of oriented planes and spheres, and leads to new remarkable incidence theorems. Most of our results are valid in hyperbolic and spherical geometries as well. We present also generalizations in spaces of higher dimension, called checkerboard IS-nets. The construction of these nets is based on a new 9 inspheres incidence theorem.

## 1 Introduction

Geometric constructions based on circles play an important role in discrete differential geometry. Circle packings and, more generally, circle patterns serve as discrete counterparts of analytic functions, see the book [20]. The origin of this idea is connected with the approach by Thurston to the Riemann mapping theorem via circle packings. Circular nets, i.e. nets with planar circular quadrilaterals, in space are discrete analogues of curvature line parametrized surfaces and orthogonal coordinate systems. They are described by discrete integrable systems, which leads to a rather developed theory, see [8]. Orthogonal circle patterns in a plane with the combinatorics of the square grid introduced by Schramm [17] are also described by discrete integrable systems. They can also be seen as special quadrilateral patterns with circumscribed quadrilaterals. Möbius geometry is a natural framework for all these theories.

In this article we construct some grids naturally related to conics and quadrics, and in particular to confocal conics.

*Institute of Science and Technology Austria (IST Austria), Am Campus 1, A – 3400 Klosterneuburg, E-mail: akopjan@gmail.com,

†Institut für Mathematik, Technische Universität Berlin, Strasse des 17. June 136, 10623 Berlin, Germany, E-mail: bobenko@math.tu-berlin.de

We consider congruences of straight lines in a plane with the combinatorics of the square grid, with all elementary quadrilaterals possessing an incircle. It follows from the Graves–Chasles theorem that the vertices of such nets (we call them incircular or IC-nets) lie on confocal conics. This gives a simple geometric construction of IC-nets starting with 2 circles and their 5 tangent lines, three of which are common (Corollary 2.10). After we presented our results at an Oberwolfach conference Serge Tabachnikov informed us about two papers of Böhm [9, 10] from the 1960s, which seem to be forgotten. In these papers Böhm gives a new proof of Ivory's theorem using the Graves–Chasles theorem. He also introduced IC-nets and the corresponding three-dimensional nets with all hypercubes possessing an insphere. Our presentation of IC-nets in Section 2 is more detailed and includes some new results.

Our main new results are on checkerboard IC-nets in the plane. These are the congruences of straight lines in the plane with the combinatorics of the square grid, combinatorially colored as a checkerboard, such that all black coordinate quadrilaterals possess inscribed circles. In Section 3 we show how this larger class of IC-nets appears quite naturally in Laguerre geometry of oriented planes. IC-nets appear then as a special case when the straight lines of the congruence coincide in pairs. The construction is based on a new incidence theorem (Theorem 3.3, see also Fig. 13) involving 13 circles and 12 straight lines. We were not able to find an elementary proof of this theorem and present a proof based on the cyclographic model of Laguerre geometry.

We also present generalizations of these nets in spaces of higher dimension, called checkerboard IS-nets. In Section 4 we show that the geometry of these nets is more rigid than the one of the planar checkerboard IC-nets. The construction is based on a new incidence theorem (Theorem 4.4) involving 9 inspheres.

Most of the results are valid in hyperbolic and spherical geometries as well. We present the corresponding modifications for the hyperbolic space in Section 5.

IC-nets are closely related to Poncelet grids (called also Poncelet-Darboux grids) introduced and studied by Darboux [11] and after him by several authors [18] and [15]. These are generated by a (Poncelet) polygon with vertices on an ellipse $\alpha$, the edges of which touch an ellipse $\alpha'$ located within $\alpha$. The straight lines of edges comprise the Poncelet grid. If the ellipses $\alpha$ and $\alpha'$ are confocal then the Poncelet grid is a periodic IC-net, see Fig. 1. On the other hand any Poncelet grid is a projective image of a periodic IC-net. The last claim follows from the fact that two nested ellipses can be made confocal by a projective transformation. As a corollary of Theorem 2.1 we obtain that the perspectivity property (vii) is valid for general Poncelet grids.

Finally, we would like to mention a number of recent attempts to discretize quadrics in general and confocal systems of quadrics in particular. In [21] a discretization of the defining property of a conic as an image of a circle under a projective transformation is considered. Since a natural discretization of a circle is a regular polygon, one ends up with a class of discrete curves that are projective images of regular polygons. Although this class has several nice geometric properties, this discretization is too simplistic.

A version of discrete confocal quadrics in any dimension was introduced in [7] in the framework of the theory of integrable systems. Starting with an integrable discretization of the Euler–Darboux system, which describes classical confocal quadrics, discrete confocal quadrics were defined analytically. They turn out to have a remarkable geometric property: all discrete two-dimensional level surfaces of the so-defined discrete quadrics are Koenigs nets. This is a very important property since together with the orthogonality condition it characterizes confocal quadrics. The corresponding discrete orthogonality condition was formulated in terms of two combinatorially dual nets.

The geometric patterns constructed in this paper are rather rigid. It would be interesting to find an appropriate analytic description in terms of difference equations. Since the geometric constructions depend on finitely many parameters, it is probably an interesting special ordinary difference equation.

**Acknowledgements.** This research was supported by the DFG Collaborative Research Center TRR 109 "Discretization in Geometry and Dynamics". The first author was also supported by People Programme (Marie Curie Actions) of the European Union's Seventh Framework Programme (FP7/2007-2013) under REA grant agreement n°[291734].

We would like to thank Serge Tabachnikov and Yuri Suris for showing us the references [9, 10, 12] where some results of Section 2 were originally obtained. We are grateful to Tim Hoffmann and Mirco Kraenz for producing Fig. 18. We also grateful to the anonymous reviewers for their careful reading of the paper and valuable suggestions.

# 2 Incircular (IC) nets

## 2.1 Main theorem

We consider maps of the square grid to the plane $f : \mathbb{Z}^2 \to \mathbb{R}^2$ and use the following notations:

- $f_{i,j} = f(i,j)$ for the vertices of the net,
- $\square^c_{i,j}$ for the quadrilateral $(f_{i,j}, f_{i+c,j}, f_{i+c,j+c}, f_{i,j+c})$ which we call a *net-square*,
- $\square_{i,j}$ for the net-square $\square^1_{i,j}$ which we call a *unit net-square.*

We denote a rectangle in $\mathbb{Z}^2$ by $\mathbb{P} = \{(i,j) \in \mathbb{Z}^2 | m_1 < i < m_2, n_1 < j < n_2\}$.

**Definition 2.1.** An IC-net (inscribed circular net) is a map $f : \mathbb{P} \to \mathbb{R}^2$ satisfying the following conditions:

(i) For any integer $i$ the points $\{f_{i,j} | j \in \mathbb{Z}\}$ lie on a straight line $\ell_i$ preserving the order, i.e the point $f_{i,j}$ lies between $f_{i,j-1}$ and $f_{i,j+1}$. The same holds for points $\{f_{i,j} | i \in \mathbb{Z}\}$ which lie on a straight line $m_j$. We call the lines $\ell_i, m_j$ the *lines of the IC-net.*

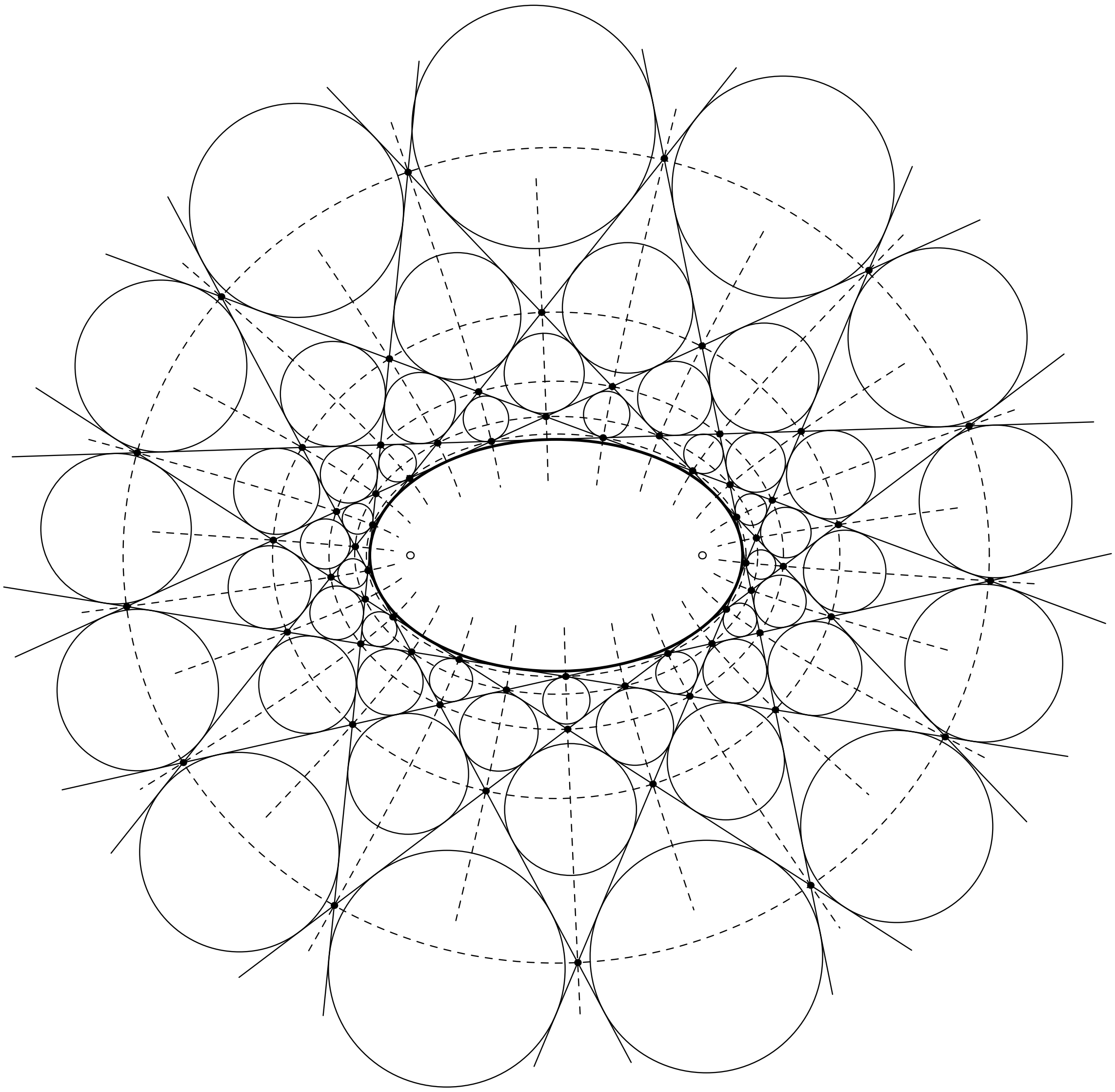

Figure 1: A Poncelet IC-net

(ii) all unit net-squares $\square_{i,j}$ are circumscribed. We denote the inscribed circle of $\square_{i,j}$ by $\omega_{i,j}$ and its center by $o_{i,j}$.

An example of an IC-net is presented in Fig. 1. IC-nets have remarkable geometric properties, which we summarize in the following theorem.

**Theorem 2.1.** *Let $f$ be an IC-net. Then the following properties hold:*

*(i) All lines of the IC-net $f$ touch some conic $\alpha$ (possibly degenerate).*

*(ii) The points $f_{i,j}$, where $i + j = \text{const}$ lie on a conic confocal with $\alpha$. As well the points $f_{i,j}$, where $i - j = \text{const}$ lie on a conic confocal with $\alpha$.*

*(iii) All net-squares of $f$ are circumscribed.*

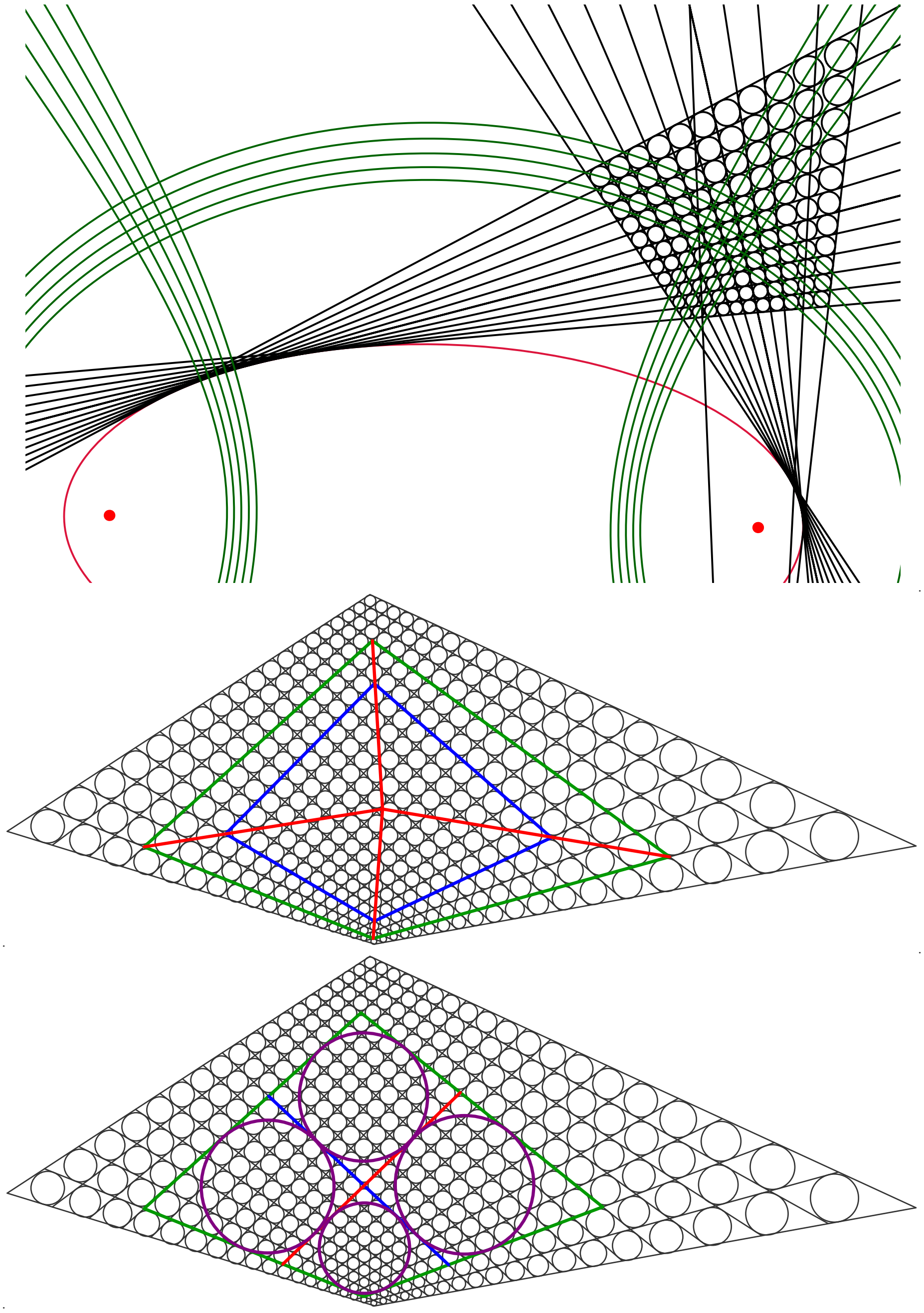

Figure 2: Geometry of IC-nets: (top) tangent lines and confocal conics, (middle) perspective net-squares, (bottom) circumscribed net-squares and equal length midlines

*(iv)* *In any net-square with even combinatorial side lengths the midlines have equal lengths:*

$$|f_{i-c,j}f_{i+c,j}| = |f_{i,j-c}f_{i,j+c}|. \tag{1}$$

*(v)* *The cross ratio*

$$cr(f_{i,j_1}, f_{i,j_2}, f_{i,j_3}, f_{i,j_4}) = \frac{(f_{i,j_1} - f_{i,j_2})(f_{i,j_3} - f_{i,j_4})}{(f_{i,j_2} - f_{i,j_3})(f_{i,j_4} - f_{i,j_1})}$$

*is independent of* $i$*. The cross ratio* $cr(f_{i_1,j}, f_{i_2,j}, f_{i_3,j}, f_{i_4,j})$ *is independent of* $j$*.*

*(vi)* *Consider the conics* $C_k$ *that contain the points* $f_{i,j}$ *with* $i + j = k$ *(see (ii)). Then for any* $l \in \mathbb{Z}$ *there exists an affine transformation* $A_{k,l} : C_k \to C_{k+2l}$ *such that* $A_{k,l}(f_{i,j}) = f_{i+l,j+l}$*. The same holds for the conics through the points* $f_{i,j}$ *with* $i - j = \mathrm{const}$*.*

*(vii)* *The net-squares* $\square^c_{i,j}$ *and* $\square^{c+2l}_{i-l,j-l}$ *are perspective.*

*(viii)* *Consider the cone in* $\mathbb{R}^3$ *intersecting the plane along the inscribed circle* $\omega_{i,j}$ *at constant oriented angle (all the apexes* $a_{i,j}$ *of these cones lie in one half-space). Then all the apexes* $a_{i,j}$ *lie on a one-sheeted hyperboloid.*

*(ix)* *All the circle centers* $o_{i,j}$ *with* $i + j = \mathrm{const}$ *lie on a conic, and* $o_{i,j}$ *with* $i - j = \mathrm{const}$ *also lie on a conic.*

*(x)* *The centers* $o_{i,j}$ *of circles of an IC-net build an affine image of an IC-net.*

We will prove this theorem in Section 2.4. But before this we present some important facts about pencils of conics.

## 2.2 The Graves–Chasles theorem

In this section we present a remarkable theorem of Graves–Chasles. It plays a crucial role in the proof of our main theorem and in construction of IC-nets. Darboux called this theorem beautiful and presented its proof in his book (see [12], p. 174). It is of course possible to prove it by direct computation (see Appendix). For completeness in this Section we give a geometric proof of the Graves–Chasles theorem and related results used for construction of IC-nets. An advantage of our approach is that it can be applied to the hyperbolic space and to the sphere as well.

**Lemma 2.2.** *Suppose on a domain* $\Omega \subset \mathbb{R}^2$ *two coordinate systems* $(x, y)$ *and* $(z, t)$ *are given. Then the following two properties are equivalent:*

*(i)* *if the points* $(x_1, y_1)$ *and* $(x_2, y_2)$ *have the same* $z$ *coordinate, then the points* $(x_1, y_2)$ *and* $(x_2, y_1)$ *have the same* $t$ *coordinate.*

*(ii)* *if the points* $(z_1, t_2)$ *and* $(z_2, t_1)$ *have the same* $x$ *coordinate, then the points* $(z_1, t_1)$ *and* $(z_2, t_2)$ *have the same* $y$ *coordinate.*

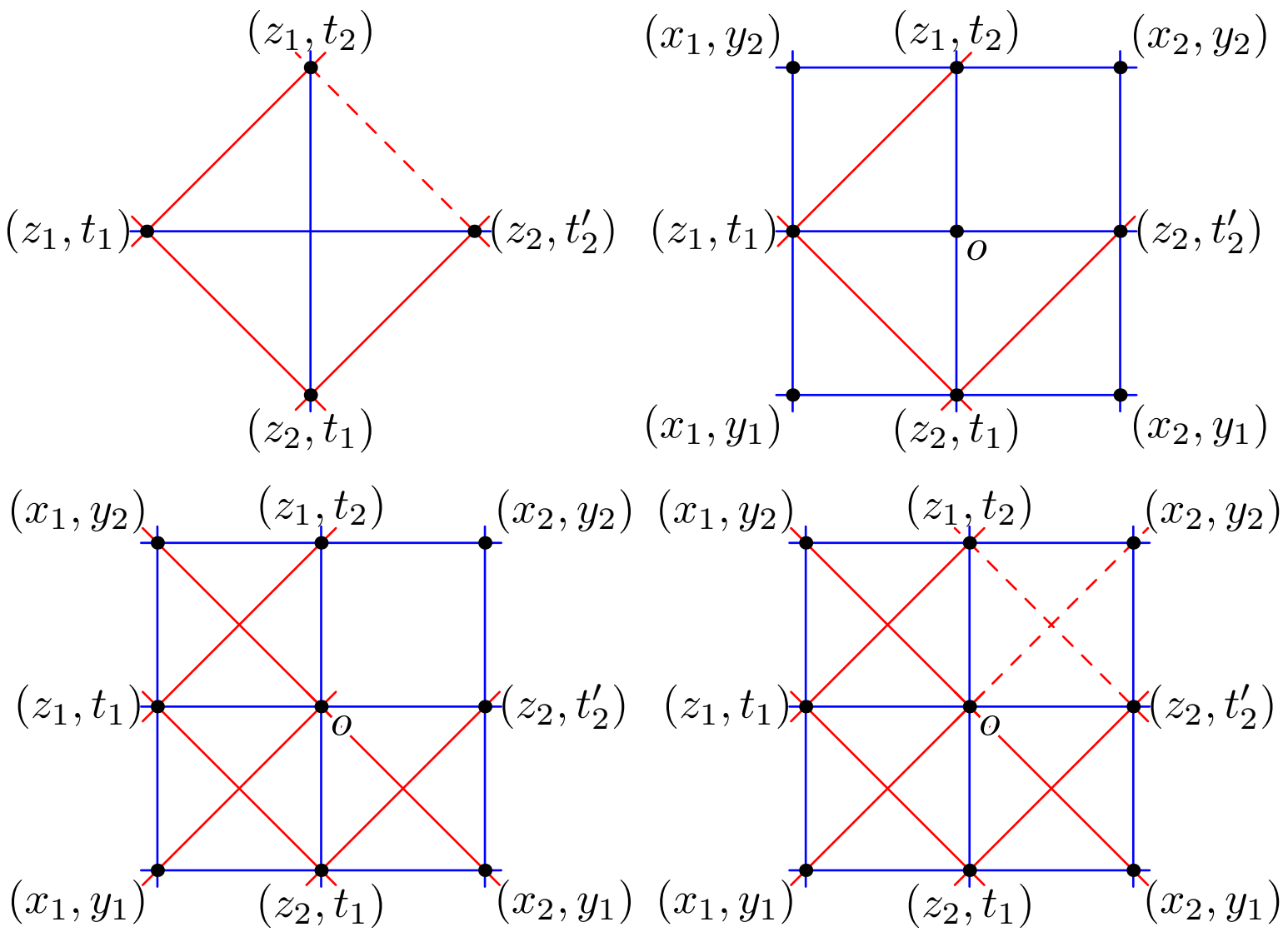

Figure 3: Two coordinate systems

*Proof.* Prove the implication (i)⇒(ii). Take two points $(z_1, t_2)$ and $(z_2, t_1)$ with the same $x$-coordinate and consider the line with the fixed $y$-coordinate passing through the point $(z_1, t_1)$ (see Fig. 3). Let $(z_2, t_2')$ be the point on this line with the coordinate $z = z_2$. We have to show that (i) implies $t_2' = t_2$.

Draw the $x, y$ coordinate lines through the points $(z_1, t_1)$, $(z_2, t_1)$, $(z_1, t_2)$, and $(z_2, t_2')$, and denote the points of their intersection as shown in Fig. 3. Applying property (i) for the pairs of points $\{(z_1, t_1), (z_1, t_2)\}$, $\{(z_2, t_1), (z_2, t_2')\}$, we see that the points $(x_1, y_2)$, $o$, and $(x_2, y_1)$ have equal $t$-coordinates. Further (i) for the pairs $\{(z_1, t_1), (z_2, t_1)\}$ and $\{(x_1, y_2), (x_2, y_1)\}$ implies that the points $(x_1, y_1)$, $o$ and $(x_2, y_2)$ have equal $z$-coordinates. Finally, applying (i) to $\{o, (x_2, y_2)\}$, we obtain that $t_2 = t_2'$.

The implication (ii)⇒(i) is proven in exactly the same way. □

**Definition 2.2.** We call a pair of coordinate line systems *diagonal-connected* if they satisfy the condition of Lemma 2.2.

We start with the following well-known lemma (for the proof see [4, Lemma 3.8] and [5, 16.6.4]).

**Lemma 2.3.** *Let $a$, $b$, $c$, $d$ be four points on a conic $\alpha$ and the lines $(ab)$ and $(cd)$ touch some other conic $\beta$. Then lines $(bd)$ and $(ac)$ touch some conic $\gamma$ from the pencil generated by the conics $\alpha$ and $\beta$. Moreover the tangent points of $\beta$ with $(ab)$ and $(cd)$, and the tangent points of $\gamma$ with $(bd)$ and $(ac)$ are collinear.*

Recall that a generic dual pencil of conics in $\mathbb{RP}^2$ is a set of conics with four common tangent lines (possibly complex). It is projectively dual to a generic pencil of conics (see [4] and [5]). The projectively dual form of Lemma 2.3 reads as follows.

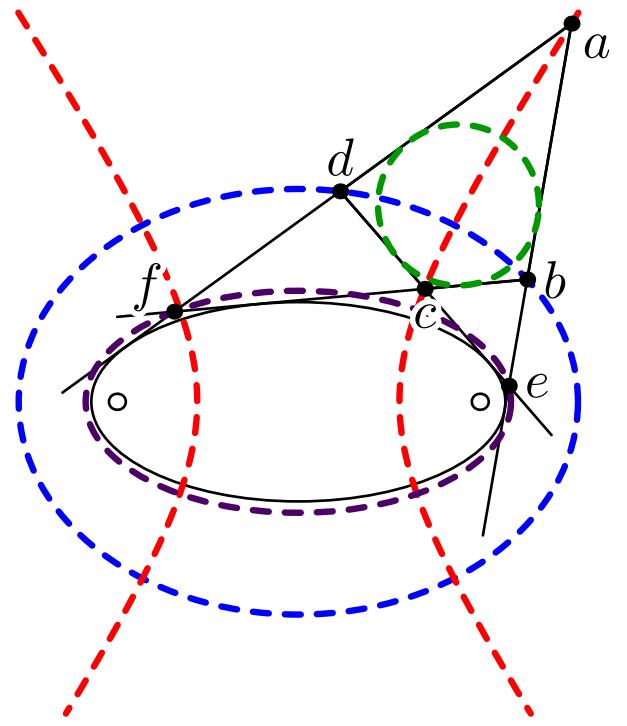

Figure 4: Graves–Chasles theorem

**Lemma 2.4.** *Consider two coordinate systems in $\Omega \subset \mathbb{R}^2$ formed by a dual pencil of conics and the family of lines tangent to some conic from this pencil respectively. These coordinate line systems are diagonal-connected.*

*Moreover, let the lines of the sides of the quadrilateral $(abcd)$ touch a conic $\alpha$ and its vertices $a$ and $c$ lie on conic $\beta$, the vertices $b$ and $d$ lie on a conic $\gamma$, and all three conics $\alpha, \beta, \gamma$ are from a dual pencil. Then the tangent lines through $a$ and $c$ to the conic $\beta$, and through $b$ and $d$ to the conic $\gamma$ intersect in one point.*

**Theorem 2.5** (Graves–Chasles theorem)**.** *Suppose that all sides of a complete quadrilateral touch a conic $\alpha$. Denote pairs of its opposite vertices by $\{a, c\}$, $\{b, d\}$, and $\{e, f\}$; see Fig. 4. Then the following four properties are equivalent:*

*(i) $(abcd)$ is circumscribed;*

*(ii) Points $a$ and $c$ lie on a conic confocal with $\alpha$;*

*(iii) Points $b$ and $d$ lie on a conic confocal with $\alpha$.*

*(iv) Points $e$ and $f$ lie on a conic confocal with $\alpha$.*

*Proof.* Note that a family of confocal conics forms a dual pencil. Therefore the equivalence of conditions (ii) and (iii) and the equivalence of conditions (ii) and (iv) follows directly from the first part of Lemma 2.4.

Let us show how the statement (i) follows from (ii) and (iii). Assume that conditions (ii) and (iii) hold. Then the second part of Lemma 2.4 implies that the tangent lines to the conics at the vertices of $(abcd)$ intersect at a point. It remains to show that these tangent lines are bisectors of the quadrilateral $(abcd)$.

The classical equal angle lemma (see, for example, [4]) states that the bisectors of the angle formed by the tangent lines from a point $p$ to the conic $\alpha$ coincide with the bisectors of the angles $\angle f_1 p f_2$, where $f_1$ and $f_2$ are the foci of the conic $\alpha$.s On the other hand the optical property of conics implies that these bisectors are tangent lines to conics confocal with $\alpha$ passing through $p$. Therefore all four tangent lines to conics at vertices of $(abcd)$ are bisectors, so their intersection point is the incenter of $(abcd)$. This completes the implication (ii)$\Rightarrow$(i).

For the proof in the opposite direction (i) $\Rightarrow$ (ii), we use the uniqueness of the configuration. Indeed, choose a point $c'$ on the line $(bc)$ such that $a$ and $c'$ lie on a conic from our confocal family, and point $d'$ on the line $(ad)$ such that $(c'd')$ is tangent to $\alpha$. We have proved already that then $(abc'd')$ is circumscribed. Moreover its incircle coincides with the incircle of the quadrilateral $(abcd)$ because it is uniquely determined by the lines $(ad)$, $(ab)$, and $(bc)$. The incircle and the conic $\alpha$ have no more than four common tangent lines. Since they already have three common tangent lines $(ad)$, $(ab)$, and $(bc)$, the common tangent lines $(cd)$ and $(c'd')$ coincide. □

We will use also the following direct corollary of this Theorem.

**Corollary 2.6.** *Let the lines $(ab)$, $(bc)$, $(cd)$ of three sides of a circumscribed quadrilateral $(abcd)$ touch a conic $\alpha$ and the vertices $a$, $c$ lie on a conic confocal to $\alpha$. Then the line $(ad)$ of the forth side also touches $\alpha$.*

*Remark* 1. Dual pencils of conics are studied in particular in [3]. From Lemma 1 of [3] it follows that if the lines of the sides of a quadrilateral $(abcd)$ touch a conic $\alpha$ and a circle with center $o$, then the lines of the sides of the quadrilateral $f_1af_2c$, where $f_1$ and $f_2$ are the foci of the conic $\alpha$, are also tangent to a circle with center $o$. However circumscribility is equivalent to the fact that $a$ and $c$ lie on a conic with foci $f_1$ and $f_2$. That shows the equivalence (i)$\Leftrightarrow$(ii).

Applying Lemma 2.4 to the previous family of confocal conics we obtain the following useful statement.

**Lemma 2.7.** *Let $\alpha_1, \alpha_2$ be two ellipses and $\gamma_1, \gamma_2$ be two hyperbolas from the same confocal family. Let $x$, $y$, $z$, $t$ be their intersection points (see Fig. 5). Then the lines $(xz)$ and $(yt)$ touch a conic confocal to $\alpha_1$, $\alpha_2$, $\gamma_1$, $\gamma_2$.*

Now the classical Ivory theorem (see, for example, [19] and [14, sect. 30.6]) follows directly from the Graves–Chasles theorem.

**Corollary 2.8** (Ivory theorem)**.** *Let $\alpha_1$, $\alpha_2$ be two ellipses and $\gamma_1$, $\gamma_2$ be two hyperbolas from the same confocal family. Let $x$, $y$, $z$, $t$ be their intersection points (see Fig. 5). Then*

$$|xz| = |yt| \tag{2}$$

*Proof.* Due to Lemma 2.7 there exists a conic $\alpha$ from the confocal family with tangent lines $(xz)$ and $(yt)$. Let us draw four more lines tangent to $\alpha$ passing through the points $x$, $y$, $z$ and $t$. They determine the intersection points $a$, $b$, $c$, $d$ as in Fig. 5. The Graves–Chasles theorem implies that four obtained quadrilaterals shown in Fig. 5 as well as the big quadrilateral $(abcd)$ are circumscribed.

Now identity (2) follows from the fact that the sum of lengths of two opposite edges of a circumscribed quadrilateral equals the semiperimeter of the quadrilateral. Using the fact that the distances between the touching points on two exterior tangent lines common to two disjoint discs are equal it is easy to see that both $2|xz|$ and $2|yt|$ are equal to the sum of semiperimeters of quadrilaterals $(xayp)$, $(ybzp)$, $(zctp)$, and $(tdxp)$ minus the semiperimeter of $(abcd)$. □

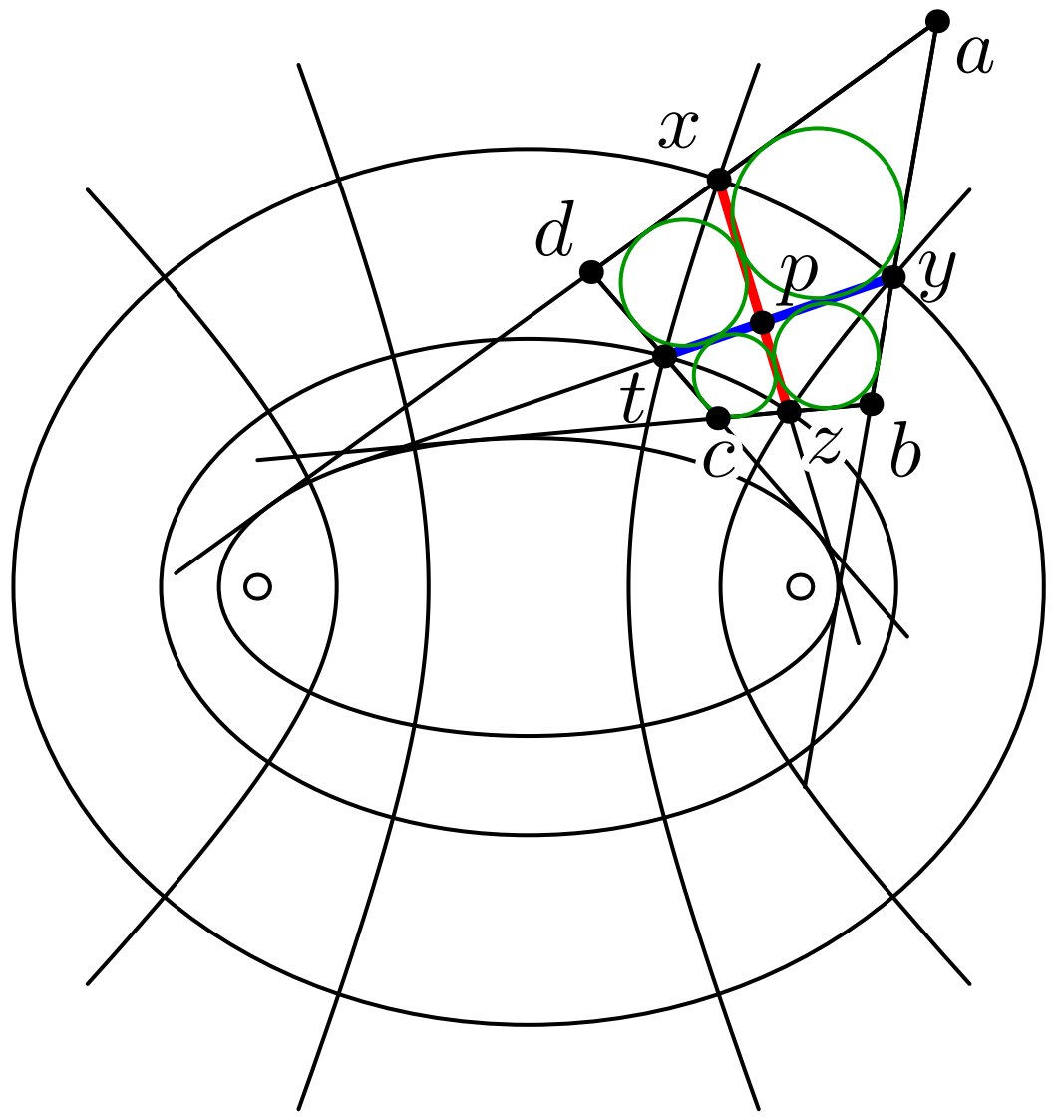

Figure 5: Proof of the Ivory theorem

We note, that the Poncelet theorem (for an exhaustive presentation see [13]), can also be proven by induction using Lemma 2.3, see [5].

We complete this section with another definition of IC-nets equivalent to Definition 2.1. Let us denote by $\ell_i$ and $m_j$ the combinatorially "vertical" and "horizontal" lines of the IC-net respectively. From Graves–Chasles theorem it follows that the intersection points $m_j \cap m_{j+1}$ and $\ell_i \cap \ell_{i+1}$ lie on the same conic from the confocal family. Varying $i$ and $j$ we obtain that this conic is independent of $i$ and $j$. Using this observation, we can define IC-net in the following way.

**Definition 2.3.** Let $\alpha$ and $\alpha'$ be confocal conics. Let $\ell_i$ and $m_j$ be lines tangent to $\alpha$ such that all the intersection points $m_j \cap m_{j+1}$ and $\ell_i \cap \ell_{i+1}$ lie on $\alpha'$ (see Fig. 8). We call $\ell_i, m_j, i, j \in \mathbb{Z}$ the lines of an IC-net and the points $f_{i,j} = \ell_i \cap m_j$ the vertices of this IC-net.

## 2.3 Construction of IC-nets

Construction of IC-nets is based on the following incidence theorem.

**Theorem 2.9** (3 × 3 incircles incidence theorem)**.** *Consider a quadrilateral which is cut in nine quadrilaterals by two pairs of lines $\ell_1, \ell_2$ and $m_1, m_2$ (see Fig. 6). Suppose all quadrilaterals except one at a corner are circumscribed. Then the ninth quadrilateral is also circumscribed.*

*Proof.* Consider the conic $\alpha$ which touches five lines $\ell_0$, $\ell_1$, $\ell_2$, $m_0$, $m_1$. Applying Corollary 2.6 several times we obtain that all lines in the figure are tangent to $\alpha$. Denote by $f_{i,j}$ the intersection points $f_{i,j} = \ell_i \cap m_j$. Theorem 2.5 implies that the pairs $f_{2,3}$ and $f_{1,2}$, $f_{1,2}$ and $f_{2,1}$, $f_{2,1}$ and $f_{3,2}$ lie on conics confocal with $\alpha$. Due to

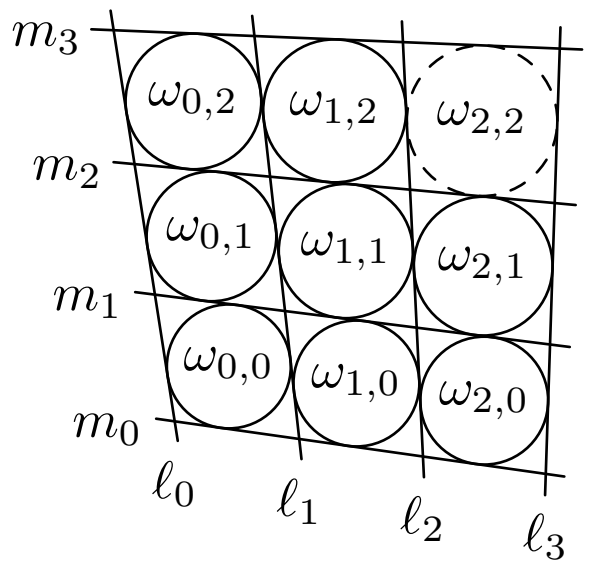

Figure 6: $3 \times 3$ incircles incidence theorem

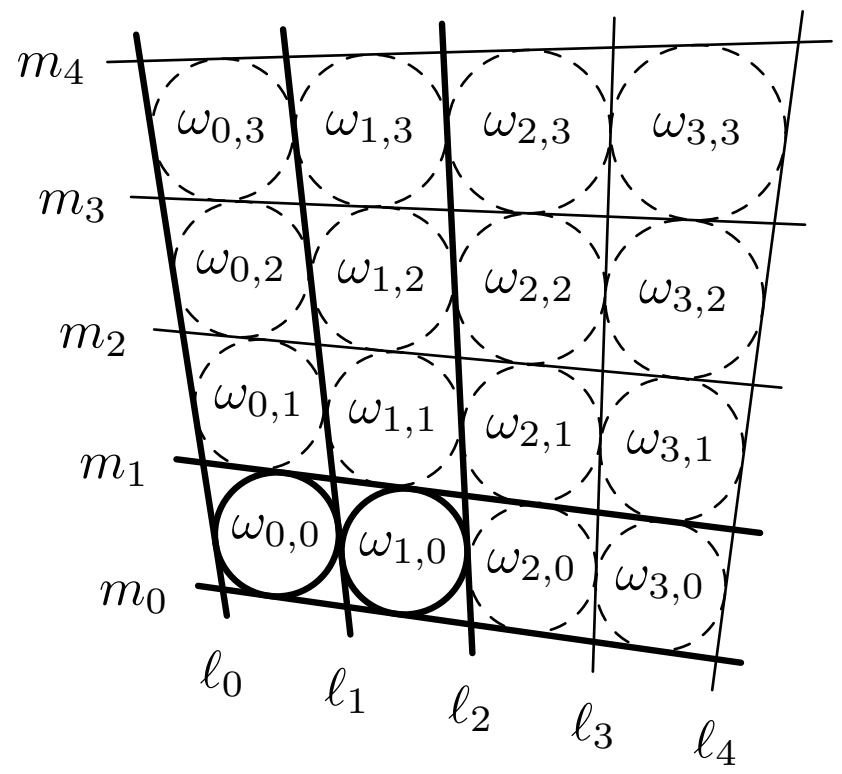

Figure 7: Construction of an IC-net from two circles and their five tangent lines

Lemma 2.7 $f_{2,3}$ and $f_{3,2}$ also lie on a conic confocal with $\alpha$. Finally Theorem 2.5 implies that the quadrilateral formed by the lines $\ell_2, \ell_3, m_2, m_3$ is circumscribed. □

There is a natural way to construct an IC-net starting from two circles.

**Corollary 2.10.** *IC-nets considered up to Euclidean motions and homothety build a real four-dimensional family. An IC-net is uniquely determined by two neighboring circles $\omega_{0,0}$, $\omega_{1,0}$ and tangent lines $\ell_0$, $\ell_1$, $\ell_2$, $m_0$, $m_1$ (see Fig. 7).*

*Proof.* Choose two non-intersecting circles $\omega_{0,0}$ and $\omega_{1,0}$ with tangent lines $\ell_0$, $\ell_1$, $\ell_2$, $m_0$, $m_1$ (see Fig. 7). Now the circles $\omega_{0,1}$, $\omega_{1,1}$, $\omega_{2,0}$ are uniquely determined. Next the common tangent line $m_2$, and further, the circles $\omega_{0,2}$, $\omega_{1,2}$, $\omega_{2,1}$ are determined. They determine the tangent lines $\ell_3$ and $m_3$. Finally the inscribed circle $\omega_{2,2}$ exists due to the incidence Theorem 2.9. Proceeding further this way one constructs the whole IC-net. □

## 2.4 Proof of Theorem 2.1

Actually we have proven already an essential part of Theorem 2.1.

(i),(ii) The corresponding properties were already proven for a $3 \times 3$ piece of an IC-net in Theorem 2.9. The global statement follows immediately.

(iii) follows directly from Theorem 2.5 and (i), (ii).

(iv) follows from the Ivory theorem (Corollary 2.8). Indeed, due to (ii) the points $f_{i-c,j}, f_{i,j-c}, f_{i+c,j}, f_{i,j+c}$ are the intersection points of two pairs of confocal ellipses and hyperbolas.

(v) It is well known (see, for example, [4]) that for any two lines tangent to a conic $\alpha$ the map from one to another generated by tangent lines to $\alpha$ is a projective map. Therefore it preserves cross-ratios of points. The map

$$(f_{i_1,j_1}, f_{i_1,j_2}, f_{i_1,j_3}, f_{i_1,j_4}) \mapsto (f_{i_2,j_1}, f_{i_2,j_2}, f_{i_2,j_3}, f_{i_2,j_4})$$

is exactly of this type.

(vi) follows from the fact that the map of a conic to another conic of the same signature from a confocal family defined through intersections by confocal conics with other signature is an affine map (see, for example, [14] or [15]). This fact follows immediately from the equations of confocal conics.

(vii) the proof is the same as the proof of property (ii) of Theorem 4.1.

(viii) follows from item (v) of Theorem 3.1.

(ix),(x) Let $\ell_i$ and $m_j$ be the combinatorially "vertical" and "horizontal" lines of an IC-net (see Definition 2.3). Let $\ell'_i$ and $m'_j$ the bisectors of the lines $\ell_i$ and $\ell_{i+1}$, and $m_j$ and $m_{j+1}$ respectively. Segments of the lines $\ell_i$ build a billiard trajectory in $\alpha'$ since they touch the confocal conic $\alpha$ (see, for example, [14] or [13]). Thus the lines $\ell'_i$ and $m'_j$ are tangent to the conic $\alpha'$.

Both claims (ix) and (x) follow from the following simple lemma.

**Lemma 2.11.** *Let $\alpha$ and $\alpha'$ be confocal conics defining an IC-net as in Definition 2.3. The lines $\ell_i$ are tangent to $\alpha$ and intersect on $\alpha'$; $\ell_i \cap \ell_{i+1} \in \alpha'$. Let $\alpha''$ be the conic that contains the intersection points $\ell'_i \cap \ell'_{i+1}$ of the bisector lines (see Fig. 8). Then there exists an affine transformation that maps $\alpha'$ to $\alpha$ and $\alpha''$ to $\alpha'$.*

*Proof.* Let

$$\frac{x^2}{a^2} + \frac{y^2}{b^2} = 1, \quad \frac{x^2}{a'^2} + \frac{y^2}{b'^2} = 1$$

be the equations of $\alpha$ and $\alpha'$ respectively. Points of $\alpha''$ are dual with respect to $\alpha'$ to lines tangent to $\alpha$. Thus, the conic $\alpha''$ is given by

$$\frac{x^2 a^2}{a'^4} + \frac{y^2 b^2}{b'^4} = 1.$$

The affine transformation $A$ that maps $\alpha'$ to $\alpha$ is given by $A(x, y) = (\lambda x, \mu y)$ with $\lambda = \frac{a}{a'}, \mu = \frac{b}{b'}$. Obviously $A$ maps $\alpha''$ to $\alpha'$. □

The affine transformation $A$ maps the centers of circles and the lines through the centers of circles of an IC-net to the vertices of an IC-net and its lines respectively. This completes the proof of the Theorem.

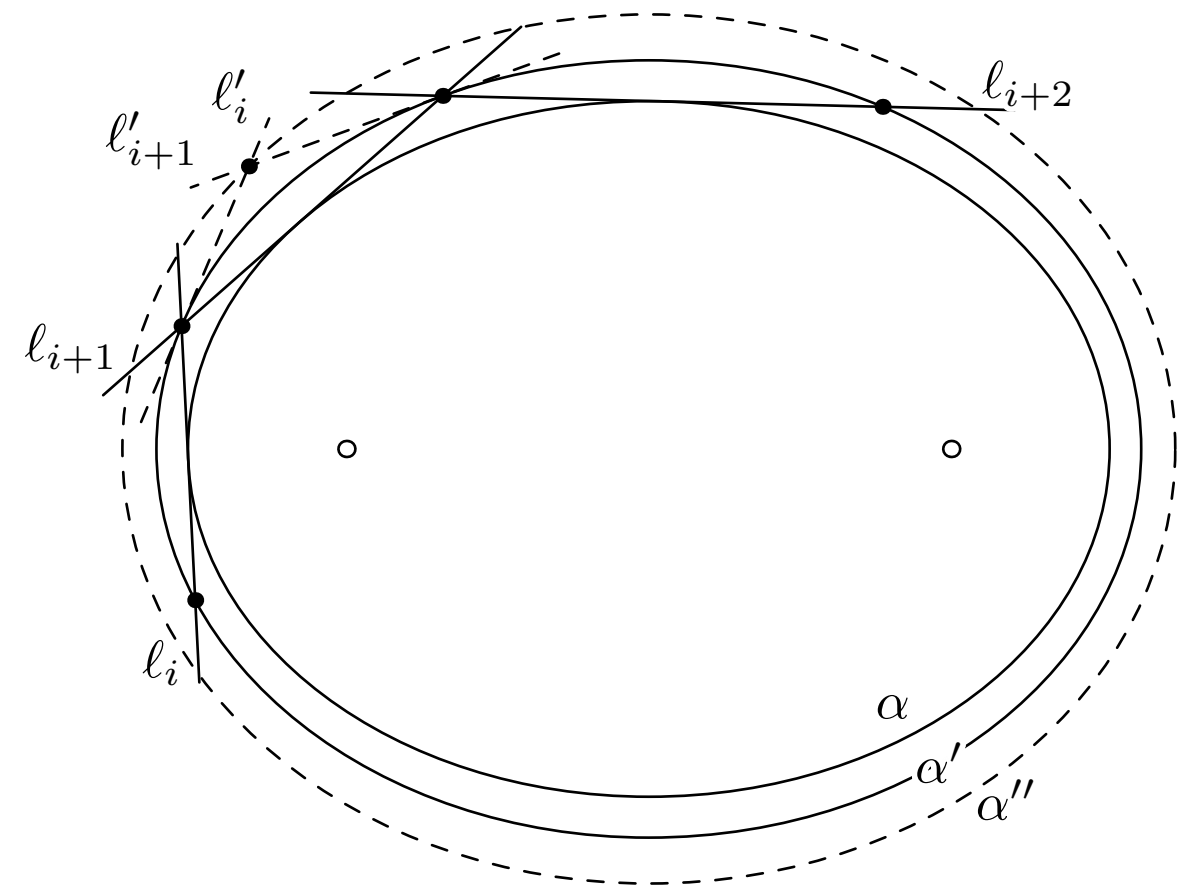

Figure 8: Three conics related by an affine transformation $A$ in Lemma 2.11: $A(\alpha'') = \alpha'$, $A(\alpha') = \alpha$. The conics $\alpha$ and $\alpha'$ are confocal.

# 3 Checkerboard IC-net

## 3.1 Definition and geometric properties of checkerboard IC-nets

**Definition 3.1.** A checkerboard IC-net is a map $f : \mathbb{P} \to \mathbb{R}^2$ satisfying the following conditions:

1. For any integer $i$ the points $\{f_{i,j} | j \in \mathbb{Z}\}$ lie on a straight line preserving the order, i.e the point $f_{i,j}$ lies between $f_{i,j-1}$ and $f_{i,j+1}$. The same holds for points $\{f_{i,j} | i \in \mathbb{Z}\}$. We call these lines the *lines of the checkerboard IC-net.*

2. For any integer $i$ and $j$ with the same parity the quadrilateral with vertices $f_{i,j}$, $f_{i+1,j}$, $f_{i+1,j+1}$, $f_{i,j+1}$ is circumscribed.

This class of nets with inscribed circles is in some sense more natural then IC-nets. The reason is that all circles and lines of a checkerboard IC-net can be consistently oriented. This shows that this class of nets belongs to Laguerre geometry, which studies oriented lines and circles that are in oriented contact (see, for example, [6]). We will use the Laguerre geometric description to prove some non-trivial incidence theorems that we have not found in the literature.

We call quadrilaterals $\square_{i,j}$ with vertices $f_{i,j}$, $f_{i+1,j}$, $f_{i+1,j+1}$, $f_{i,j+1}$ with even $i+j$ *unit net-squares* of checkerboard IC-net. The quadrilaterals $\square_{i,j}^c$ with vertices $f_{i,j}$, $f_{i+c,j}$, $f_{i+c,j+c}$, $f_{i,j+c}$ with even $i+j$ and odd $c$ we call *net-squares.*

**Theorem 3.1.** *Let $f$ be a checkerboard IC-net. Then the following properties hold:*

*(i) All net-squares are circumscribed (Fig. 9).*

*(ii) Net-squares $\square_{i,j}^c$ and $\square_{i-l,j-l}^{c+2l}$, where $l$ is odd, are perspective (Fig. 10).*

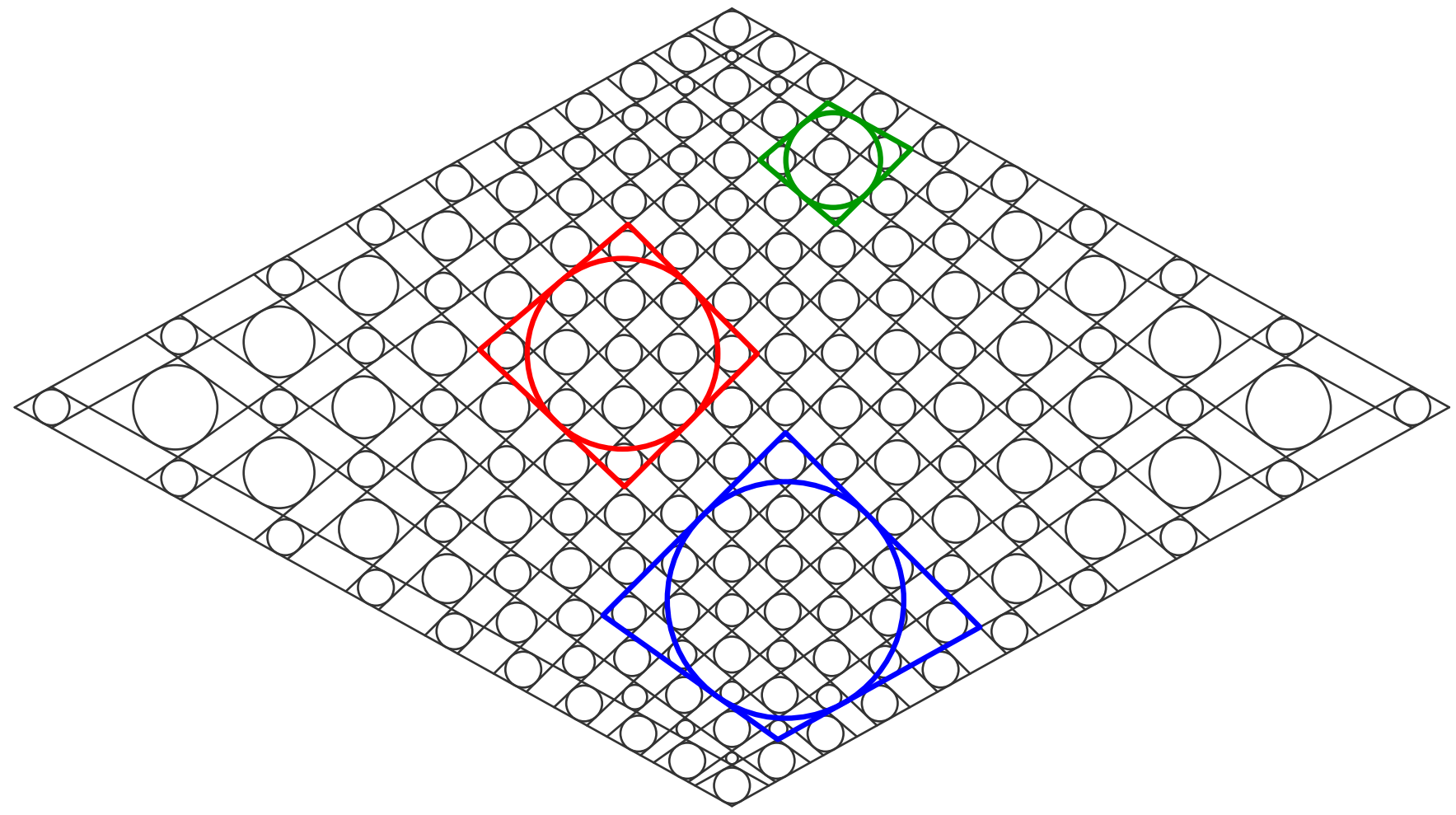

Figure 9: Checkerboard IC-nets: circumscribed net-squares

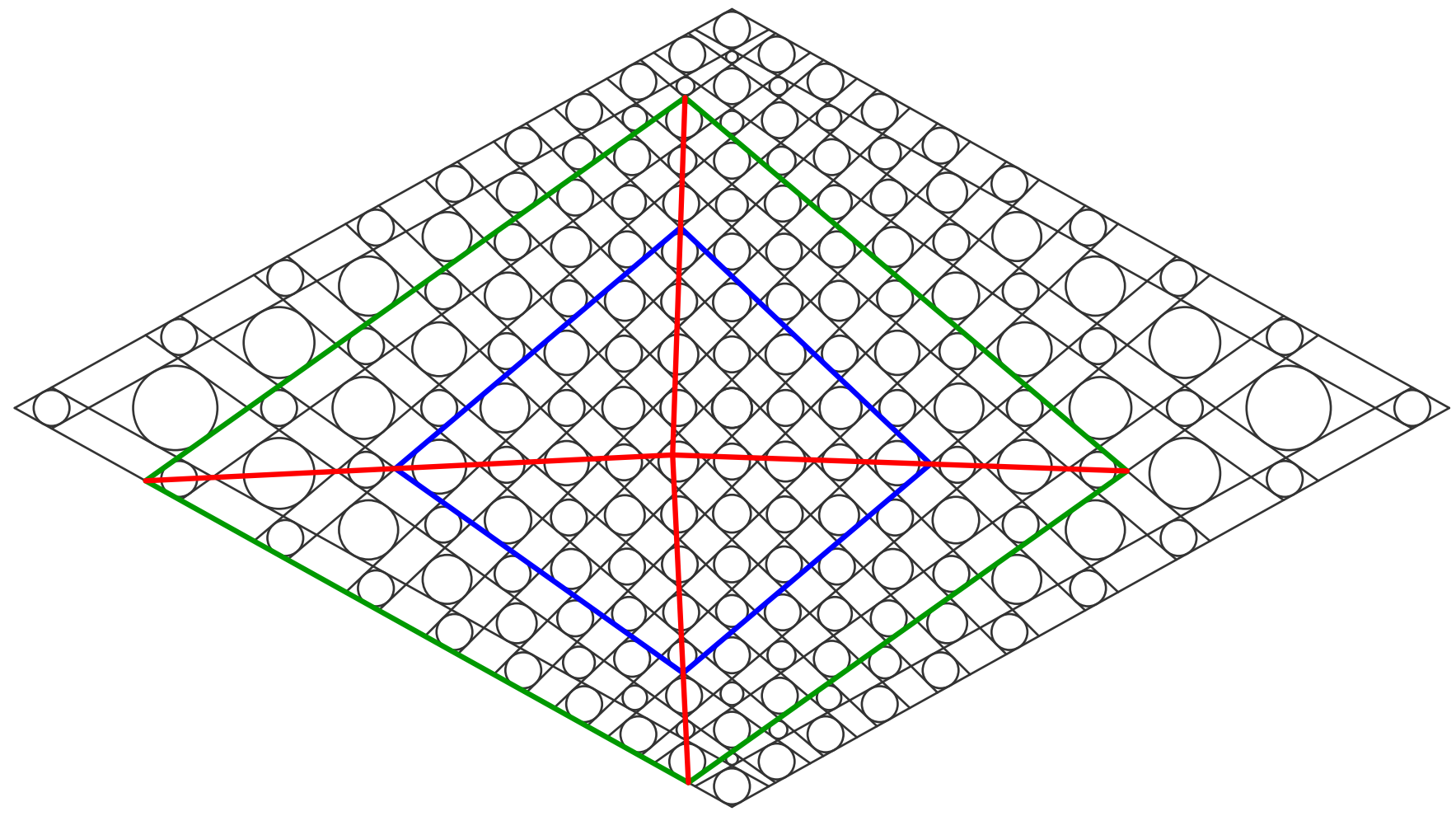

Figure 10: Checkerboard IC-nets: perspective net-squares

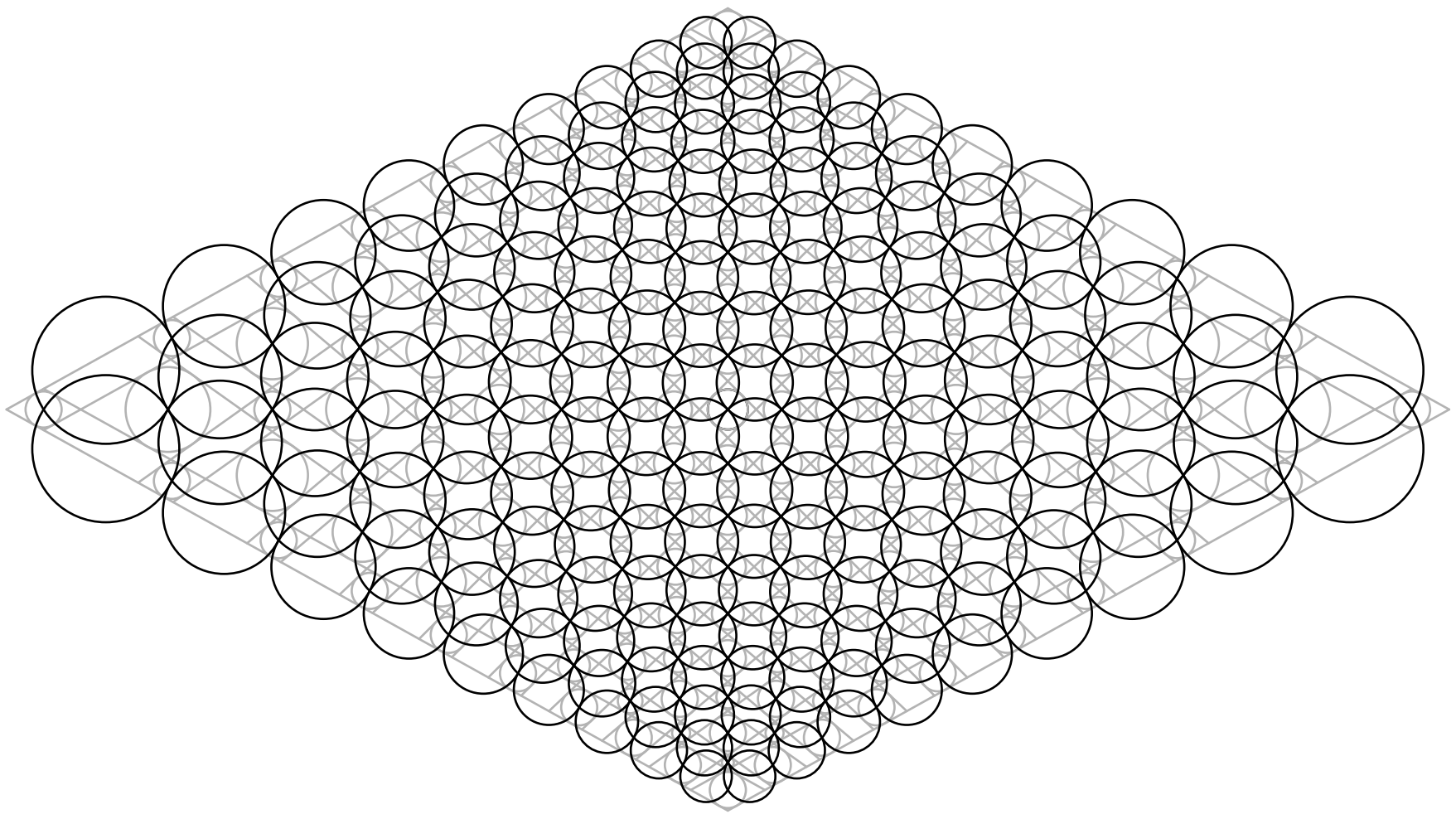

Figure 11: Checkerboard IC-nets as circular-conical nets

*(iii)* *The points $f_{i,j}$, where $i+j$ is an odd constant lie on a conic. The points $f_{i,j}$, where $i-j$ is an even constant lie on a conic as well.*

*(iv)* *(Ivory-type theorem) We define the distance $d_C(\Box_{a,b}, \Box_{c,d})$ between two unit net-squares $\Box_{a,b}$ and $\Box_{c,d}$ of a checkerboard net as the distance between the tangent points on a common exterior tangent line to the circles $\omega_{a,b}$ and $\omega_{c,d}$ inscribed in $\Box_{a,b}$ and $\Box_{c,d}$ respectively. In case $a=c$ or $b=d$ these tangent lines are the lines of the checkerboard IC-net.*

*For any $(i,j)\in\mathbb{Z}^2$, with even $i+j$ and any integer even $c$ one has*

$$d_C(\Box_{i-c,j}, \Box_{i+c,j}) = d_C(\Box_{i,j-c}, \Box_{i,j+c}). \tag{3}$$

*(v)* *Let $\omega_{i,j}$ be the inscribed circle of the unit net-square $\Box_{i,j}$. Consider the cone in $\mathbb{R}^3$ intersecting the plane along $\omega_{i,j}$ at constant oriented angle (all the apexes $a_{i,j}$ of these cones lie in one half-space). Then all the apexes $\{a_{i,j}|i+j=4n, n\in\mathbb{Z}\}$ lie on one-sheeted hyperboloid. The the apexes $\{a_{i,j}|i+j=4n+2, n\in\mathbb{Z}\}$ lie on one-sheeted hyperboloid as well.*

*(vi)* *The centers $o_{i,j}$ of the incircles of a checkerboard IC-net build a circle-conical net, i.e. a net that is simultaneously circular and conical (see [8]). Recall that circular nets are the nets with circular quadrilaterals $(o_{i,j}o_{i+1,j+1}o_{i,j+2}o_{i-1,j+1})$ and conical nets in plane are characterized by the condition that the sums of two opposite angles at a vertex are equal (and equal to $\pi$).*

We start with the following classical Lemma which can be found for example in [23, Sec. 67].

**Lemma 3.2.** *Consider a quadrilateral that is cut in nine quadrilaterals by two pairs of lines (see Fig. 12). Suppose the center quadrilateral and all the corner quadrilaterals are circumscribed. Then the "big" quadrilateral is also circumscribed.*

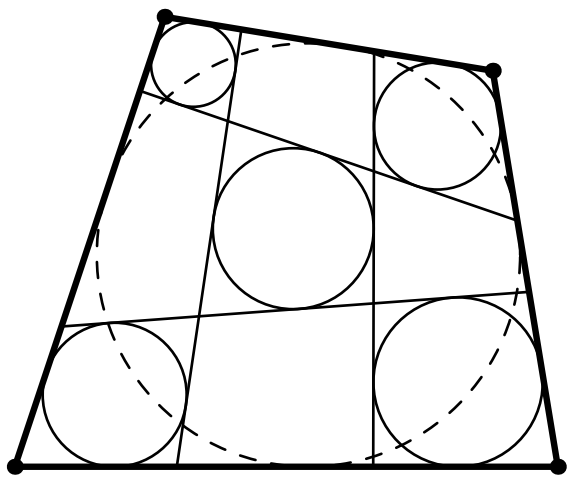

Figure 12: Six incirlces lemma

An elementary proof of this lemma is based on the fact that the distances between the touching points on two exterior tangent lines common to two disjoint discs are equal. Using this fact one can show that the differences of the sums of the lengths of the opposite sides of the central "small" and "big" quadrilaterals are equal. This implies that they are simultaneously circumscribed.

*Proof.* of Theorem 3.1

(i) We will prove this by induction on $c$. From the definition we obtain the claim for $c = 1$. Suppose the circumscribility is known for all net-squares of size $c$. Let us prove it for net-squares of size $c + 2$. Applying Lemma 3.2 for the unit net-squares $\square_{i,j}$, $\square_{i,j+c+1}$, $\square_{i+c+1,j+c+1}$, $\square_{i+c+1,j}$ and the net-square $\square^c_{i+1,j+1}$ we get that the net-square $\square^{c+2}_{i,j}$ is also circumscribed.

(ii) The proof is the same as the proof of property (ii) of Theorem 4.1.

(iii) We will prove the claim for vertices $f_{i,j}$ with $i - j = 0$, for other cases the proof is the same. For that we will show that any six successive points lie on a some conic. Without loss of generality, we assume that $i = 1$. By the Pascal theorem it is enough to show that the following three points of intersection of lines $(f_{0,0}f_{1,1}) \cap (f_{3,3}f_{4,4})$, $(f_{1,1}f_{2,2}) \cap (f_{4,4}f_{5,5})$, $(f_{2,2}f_{3,3}) \cap (f_{5,5}f_{0,0})$ are collinear.

The point $f_{0,0}$ is the center of positive homothety of incircles of the net-squares $\square_{0,0}$ and $\square^3_{0,0}$. The point $f_{1,1}$ is a center of negative homothety of incircles of the net-squares $\square_{0,0}$ and $\square^3_{1,1}$. Therefore, by the Monge theorem (see [22, Sec. 20]), the line $(f_{0,0}f_{1,1})$ passes through the center of negative homothety of incircles of the net-squares $\square^3_{1,1}$ and $\square^3_{0,0}$. Analogously the line $(f_{3,3}f_{4,4})$ passed through this center. We obtain that the point $(f_{0,0}f_{1,1}) \cap (f_{3,3}f_{4,4})$ is the center of negative homothety of the incircles of $\square^3_{1,1}$ and $\square^3_{0,0}$.

Using the same argument we can prove that $f_{1,1}f_{2,2} \cap f_{4,4}f_{5,5}$ is the center of negative homothety of incircles of $\square^3_{1,1}$ and $\square^3_{2,2}$. The point $(f_{2,2}f_{3,3}) \cap (f_{5,5}f_{0,0})$ is the center of positive homothety of incircles of $\square^3_{0,0}$ and $\square^3_{2,2}$. Applying the Monge theorem again we obtain that these three centers of homotheties lie on a one line.

(iv) Since the net-square $\square^{c+1}_{i,j}$ is circumscribed the sums of the lengths of its opposite sides are equal. The sum of the lengths of two opposite sides of $\square^{c+1}_{i,j}$ is equal to $d_C(\Box_{i,j}, \Box_{i+c,j}) + d_C(\Box_{i,j+c}, \Box_{i+c,j+c})$ plus the sum of the lengths of the intervals from the corners of $\Box_{i,j}^{c+1}$ to the touching points with the inscribed circles of corresponding corner unit net-squares. We obtain

$$d_C(\Box_{i,j}, \Box_{i+c,j}) + d_C(\Box_{i,j+c}, \Box_{i+c,j+c}) = d_C(\Box_{i,j}, \Box_{i,j+c}) + d_C(\Box_{i+c,j}, \Box_{i+c,j+c}).$$

Applying this equality to $\Box_{i-c,j-c}^{c+1}$, $\Box_{i-c,j}^{c+1}$, $\Box_{i,j-c}^{c+1}$, $\Box_{i,j}^{c+1}$ and $\Box_{i-c,j-c}^{2c+1}$ we get

$$\begin{aligned}
2d_C(\Box_{i-c,j}, \Box_{i+c,j}) &= 2d_C(\Box_{i-c,j}, \Box_{i,j}) + 2d_C(\Box_{i,j}, \Box_{i+c,j}) = \\
&= d_C(\Box_{i-c,j-c}, \Box_{i-c,j}) + d_C(\Box_{i,j-c}, \Box_{i,j}) - d_C(\Box_{i-c,j-c}, \Box_{i,j-c}) + \\
&+ d_C(\Box_{i-c,j}, \Box_{i-c,j+c}) + d_C(\Box_{i,j}, \Box_{i,j+c}) - d_C(\Box_{i-c,j}, \Box_{i,j}) + \\
&+ d_C(\Box_{i,j-c}, \Box_{i,j}) + d_C(\Box_{i+c,j-c}, \Box_{i+c,j}) - d_C(\Box_{i,j-c}, \Box_{i+c,j-c}) + \\
&+ d_C(\Box_{i,j}, \Box_{i,j+c}) + d_C(\Box_{i+c,j}, \Box_{i+c,j+c}) - d_C(\Box_{i,j}, \Box_{i+c,j}) = \\
&= d_C(\Box_{i-c,j-c}, \Box_{i-c,j+c}) + d_C(\Box_{i+c,j-c}, \Box_{i+c,j+c}) - \\
&- d_C(\Box_{i-c,j-c}, \Box_{i+c,j-c}) - d_C(\Box_{i-c,j+c}, \Box_{i+c,j+c}) + 2d_C(\Box_{i,j-c}, \Box_{i,j+c}) = \\
&= 2d_C(\Box_{i,j-c}, \Box_{i,j+c}). \quad (4)
\end{aligned}$$

(v) The apexes of the cones of unit net-squares with fixed $i$ (or fixed $j$) are collinear. The lines of different types (with fixed $i$ or fixed $j$) intersect each other. Therefore they are asymptotic lines (of two different families) of a one-sheeted hyperboloid.

(vi) Both angle conditions of circularity and of conicality follow immediately.

□

## 3.2 Construction of checkerboard IC-nets

The following theorem is a generalization of Theorem 2.9.

**Theorem 3.3** (checkerboard incircles incidence theorem)**.** *Consider a quadrilateral that is cut by two sets of four lines in* 25 *quadrilaterals. Color the quadrilaterals in a checkerboard pattern with black quadrilaterals at the corners. Assume that all black quadrilaterals except one at a corner are circumscribed. Then the last black quadrilateral at the corner (thirteenth quadrilateral) is also circumscribed (Figure 13).*

Before we prove this theorem let us make an important comment. As we have already pointed out in Section 3.1 checkerboard IC-nets can be oriented in such a way that their circles and lines are in oriented contact. They can be naturally described in frames of Laguerre geometry. Let us briefly introduce the cyclographic model of Laguerre geometry in the plane (see [6]).

In this model the space of oriented circles $C = \{x \in \mathbb{R}^2 | |x - c|^2 = r^2\}$ is in one-to-one correspondence with the points $a = (c, r)$ of the Minkowski space $\mathbb{R}^{2,1}$. They can be seen as the apexes of the cones of revolution intersecting the plane $\mathbb{R}^2 \subset \mathbb{R}^{2,1}$ at the angle $\pi/4$ along the circles $C$. The oriented lines $\ell \in \mathbb{R}^2$ are modelled as oriented planes $L \subset \mathbb{R}^{2,1}$ intersecting the plane $\mathbb{R}^2$ along the lines $\ell$ at the angle $\pi/4$.

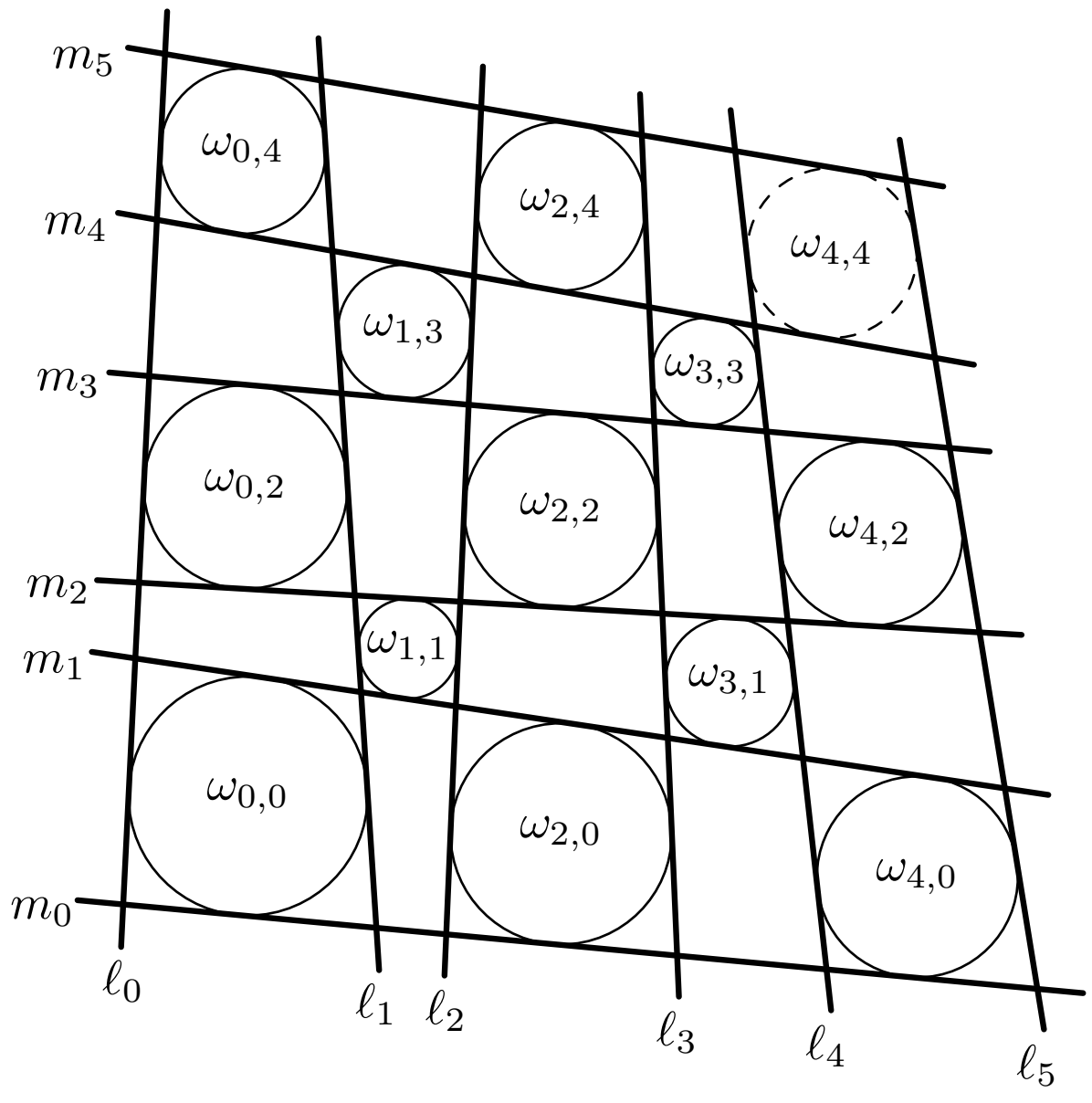

Figure 13: Checkerboard incircles incidence theorem

An oriented circle $C$ is in oriented contact with an oriented line $\ell$ if and only if $a \in L$. Two oriented circles $C_1$, $C_2 \subset \mathbb{R}^2$ are in oriented contact if and only if their representatives in the Minkowski space $a_1$, $a_2 \in \mathbb{R}^{2,1}$ differ by an isotropic vector $|a_1 - a_2| = 0$.

In this model a one-parameter family of circles that are in oriented contact to two oriented lines is represented by a straight line in $\mathbb{R}^{2,1}$. The points of this line are the apexes of the corresponding cones.

*Proof.* Let us orient the circles $\omega_{i,j}$ with even $i$, $j$ positively, the circles $\omega_{i,j}$ with odd $i$, $j$ negatively, and the tangent lines $\ell, m$ so that they are in oriented contact to these circles. Due to Lemma 3.2 the net-square $\square^3_{1,1}$ is also circumscribed, denote its inscribed circle as $\omega^3_{1,1}$ and orient it so that it is in oriented contact to the lines.

Consider the Laguerre geometry of this pattern in the cyclographic model described above. Let $a_{i,j}$ be the apex of the cone which intersect the plane at the angle $\pi/4$ along $\omega_{i,j}$. Denote by $a'_{2,2}$ the apex on the cone corresponding to the circle $\omega^3_{1,1}$. The orientation described above implies that the points $a_{1,1}$, $a_{1,3}$, $a_{3,1}$, $a_{3,3}$, $a'_{2,2}$ lie in one halfspace of $\mathbb{R}^{2,1}$ (let us say, positive third component $r$), and the points $a_{i,j}$ with even $i, j$ in the other halfspace (negative $r$).

Apexes $a$'s are collinear if and only if the corresponding circles share two common tangent lines. The following triples of points are collinear $\{a'_{2,2}, a_{1,1}, a_{0,0}\}$, $\{a'_{2,2}, a_{1,3}, a_{0,4}\}$, $\{a'_{2,2}, a_{3,1}, a_{4,0}\}$, $\{a_{2,0}, a_{2,2}, a_{2,4}\}$, and $\{a_{0,2}, a_{2,2}, a_{4,2}\}$. Moreover if the circles share a common tangent line then the corresponding apexes are coplanar. And the following quintuples of points are coplanar $\{a_{0,2}, a_{2,2}, a_{4,2}, a_{1,1}, a_{3,1}\}$, $\{a_{0,2}, a_{2,2}, a_{4,2}, a_{1,3}, a_{3,3}\}$, $\{a_{2,0}, a_{2,2}, a_{2,4}, a_{1,1}, a_{1,3}\}$, $\{a_{2,0}, a_{2,2}, a_{2,4}, a_{3,1}, a_{3,3}\}$.

We obtain an octahedron[1] $a_{2,2}, a_{1,1}, a_{3,1}, a_{3,3}, a_{1,3}, a'_{2,2}$ as in Fig. 14. The inter-section line of the face planes $(a_{1,1}a_{1,3}a_{2,2})$ and $(a_{3,1}a_{3,3}a_{2,2})$ intersects the planes $(a_{1,1}a_{3,1}a'_{2,2})$ and $(a_{3,3}a_{1,3}a'_{2,2})$ at the points $a_{2,0}$ and $a_{2,4}$ respectively. Analogously,

$$\begin{aligned} a_{0,2} &= (a_{1,1}a_{3,1}a_{2,2}) \cap (a_{1,3}a_{3,3}a_{2,2}) \cap (a_{1,1}a_{1,3}a'_{2,2}) \\ a_{4,2} &= (a_{1,1}a_{3,1}a_{2,2}) \cap (a_{1,3}a_{3,3}a_{2,2}) \cap (a_{3,1}a_{3,3}a'_{2,2}) \end{aligned} \tag{5}$$

[1]Here by octahedron we mean a polytope with combinatorics of the regular octahedron.

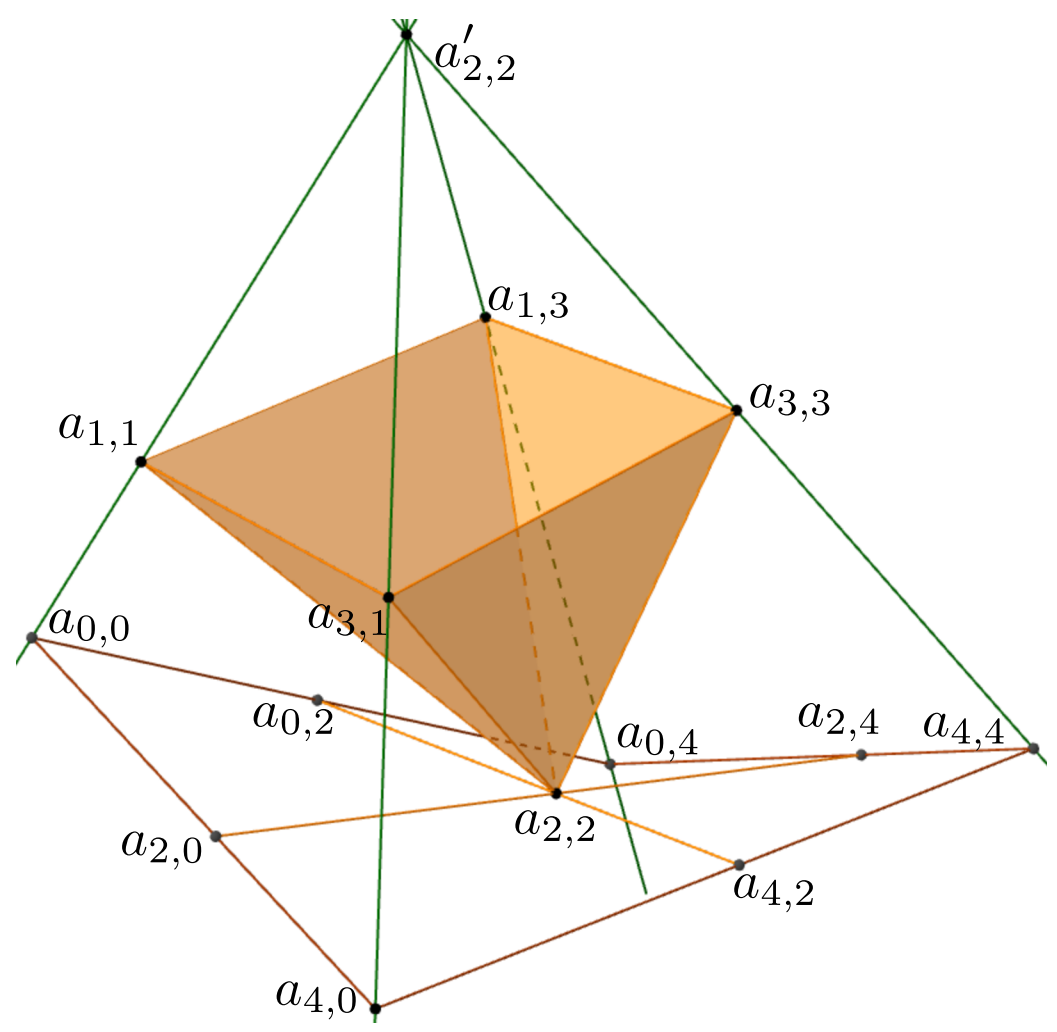

Figure 14: Projective octahedron incidence theorem

The rest of the proof follows from Theorem 3.4. □

**Theorem 3.4** (Octahedron incidence theorem)**.** *Consider an octahedron as in Fig. 14 and define the points $a_{0,2}$, $a_{2,0}$, $a_{4,2}$, $a_{2,4}$ as the intersection points of the corresponding face planes like (5). Choose an arbitrary point $a_{0,0} \in (a_{1,1}a'_{2,2})$, determine $a_{4,0} := (a_{0,0}a_{2,0}) \cap (a_{3,1}a'_{2,2})$ and $a_{0,4} := (a_{0,0}a_{0,2}) \cap (a_{1,3}a'_{2,2})$. Then the lines $(a_{0,4}a_{2,4})$ and $(a_{4,0}a_{4,2})$ intersect the line $(a_{3,3}a'_{2,2})$ in a common point (which we denote by $a_{4,4}$).*

*Proof.* Since the claim is projective we may simplify the presentation by mapping the lines $(a_{2,0}a_{2,4})$ and $(a_{0,2}a_{4,2})$ to the infinity plane. Then point $a_{2,2}$ also lies in this plane. The planes $(a_{1,1}a_{3,1}a_{2,2})$ and $(a_{1,3}a_{3,3}a_{2,2})$ become two parallel planes, the planes $(a_{1,1}a_{1,3}a_{2,2})$ and $(a_{3,1}a_{3,3}a_{2,2})$ are parallel as well.

The straight line $(a_{0,0}a_{4,0})$ is the intersection line of the planes $(a_{1,1}a_{3,1}a'_{2,2})$ and $(a_{0,0}a_{2,0}a_{2,2})$. After our normalization the latter becomes the plane parallel to $(a_{1,1}a_{1,3}a_{2,2})$ (and $(a_{3,1}a_{3,3}a_{2,2})$) and passing through $a_{0,0}$. Finally we obtain the lines $(a_{0,0}a_{4,0})$, $(a_{0,4}a_{2,4})$ and $(a_{0,0}a_{0,4})$, $(a_{4,0}a_{4,2})$ as the intersections of the corresponding face planes at $a'_{2,2}$ with the planes parallel to $(a_{1,1}a_{1,3}a_{2,2})$ and $(a_{1,1}a_{3,1}a_{2,2})$ respectively.

Now projecting the whole geometry to a plane (transversal to the line $(a_{2,2}a_{1,1})$ etc.) we obtain the incidence statement from Lemma 3.5 in plane geometry. □

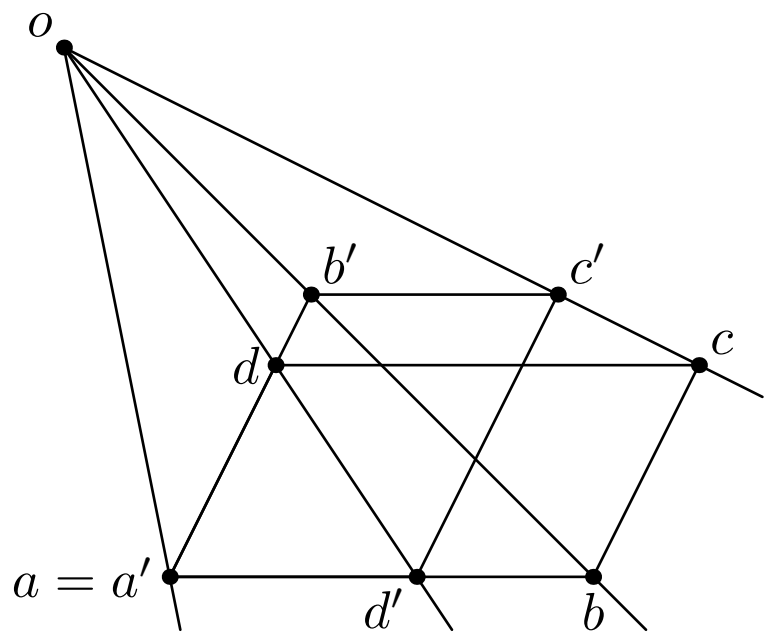

Figure 15: Proof of Lemma 3.5 from Pappus theorem

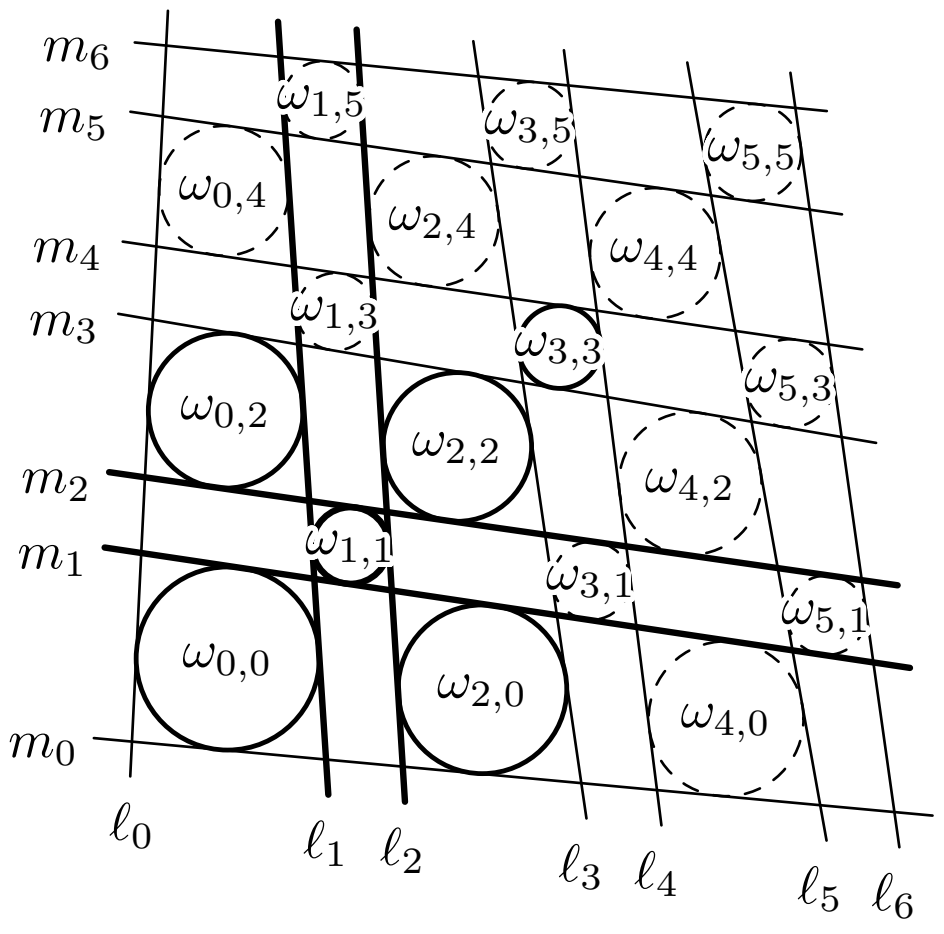

Figure 16: Construction of a checkerboard IC-net

**Lemma 3.5.** *Let $(abcd)$ be a parallelogram and $o$ be a point in plane that does not lie on the lines of the sides of the parallelogram. Let $a'$ be a point on $(oa)$. Then there exists a unique parallelogram $(a'b'c'd')$ such that lines $(bb')$, $(cc')$, $(dd')$ pass through the point $o$ and the non-corresponding sides are parallel to the sides of the original parallelogram: $(a'd') \parallel (ab)$ and $(c'd') \parallel (bc)$.*

*Proof.* Without loss of generality one can assume $a' = a$. After that the claim is just the Pappus theorem for points $b$, $c$, $d$, $b'$, $c'$, $d'$, $o$ and two points at infinity. □

The dual version of this Lemma formulated for conics can be found in [1].

**Corollary 3.6.** *Checkerboard IC-nets considered up to Euclidean motions and homothety build a real eight-dimensional family. A checkerboard IC-net is uniquely determined by five neighboring circles $\omega_{0,0}$, $\omega_{2,0}$, $\omega_{0,2}$, $\omega_{2,2}$, $\omega_{1,1}$ and circle $\omega_{3,3}$ (see Fig. 16).*

*Proof.* For constructing of a checkerboard IC-net we start with a circle $\omega_{1,1}$ and its four tangents $\ell_1$, $\ell_2$, $m_1$ and $m_2$. Then we inscribe circles $\omega_{0,0}$, $\omega_{0,2}$, $\omega_{2,2}$ and $\omega_{2,0}$ in the corners (see Figure 16). The four common tangents lines $\ell_0$, $\ell_3$, $m_0$, and $m_3$ are uniquely determined. Further, the circles $\omega_{1,3}$ and $\omega_{3,1}$ are uniquely determined.

By choosing the circle $\omega_{3,3}$ we have one degree of freedom. Now, the whole net is fixed. Indeed, we consequently determine the lines $\ell_4$ and $m_4$, then the circles $\omega_{4,0}$, $\omega_{4,2}$, $\omega_{0,4}$, $\omega_{2,4}$ and lines $\ell_5$ and $m_5$. The existence of the circle $\omega_{4,4}$ follows from the incidence theorem 3.3. Applying this theorem again and again one generates the whole checkerboard IC-net. □

## 3.3 Confocal checkerboard IC-net

IC-nets can be considered as special checkerboard IC-nets. Indeed, if for every second line and every second column of a checkerboard IC-net the incircles degenerate to points, then the lines of the net merge in pairs, all non-circumscribed quadrilaterals disappear, and one obtains an IC-net.

There is however an interesting class of IC-nets, which lies between the two we have considered. These are special checkerboard IC-nets related to confocal conics.

**Definition 3.2.** We call a checkerboard IC-net *confocal* if all lines of it are tangent to a conic.

This class is a natural generalization of IC-nets introduced in Section 2. However in contrast to IC-nets here all circles and lines can be oriented so that and the corresponding circles and lines are in oriented contact. This class can be studied in Laguerre geometry.

The following important geometric property of confocal checkerboard IC-nets can be proven exactly in the same way as the corresponding theorems in Sections 2.1 and 3.1.

**Theorem 3.7.** *Let $f$ be a confocal checkerboard IC-net all lines of which are tangent to a conic $\alpha$. Then the points $f_{i,j}$, where $i + j$ is an odd constant lie on a conic confocal to $\alpha$. The points $f_{i,j}$, where $i - j$ is an even constant lie on a conic confocal to $\alpha$ as well.*

This allows us to define confocal checkerboard IC-nets through confocal conics similarly to Definition 2.3 of IC-nets.

**Definition 3.3.** Let $\alpha$, $\alpha'$ and $\alpha''$ be confocal conics. Let $\ell_i$ and $m_i$ be lines tangent to $\alpha$ such that all the points $\ell_i \cap \ell_{i+1}$ and $m_i \cap m_{i+1}$ lie on $\alpha'$ for odd $i$ and on $\alpha''$ for even $i$. A confocal checkerboard IC-net is a map $f : \mathbb{Z}^2 \to \mathbb{R}^3$ given by the intersection points $f_{i,j} = \ell_i \cap m_j$.

# 4 Checkerboard inspherical nets in $\mathbb{R}^3$

## 4.1 Definition and geometric properties of checkerboard IS-nets

In this section we introduce a natural three-dimensional version of checkerboard IC-nets. All results and proofs can be generalized for higher dimensions. We consider images of the integer lattice $f : \mathbb{Z}^3 \to \mathbb{R}^3$ or of a cuboid $\mathbb{P} \subset \mathbb{Z}^3$ of full dimension. Let us denote $\square^c_{i,j,k}$ the cube with the vertices $f_{i,j,k}$, $f_{i+c,j,k}$, $f_{i,j+c,k}$, $f_{i,j,k+c}$,$f_{i+c,j+c,k}$, $f_{i,j+c,k+c}$, $f_{i+c,j,k+c}$, $f_{i+c,j+c,k+c}$.

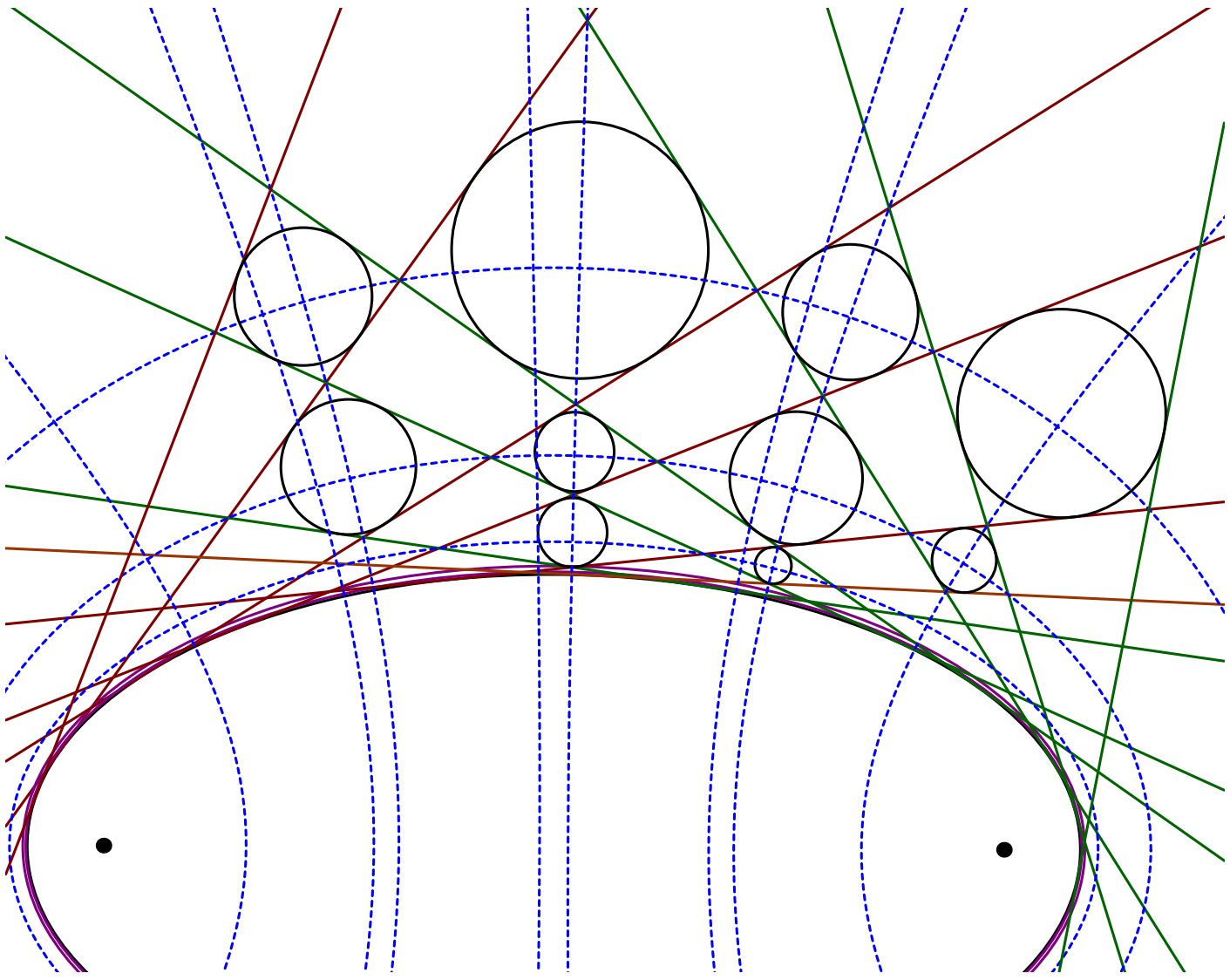

Figure 17: To a definition of confocal checkerboard IC-nets

We call $\square^c_{i,j,k}$ a *net cube* if $i, j, k$ are all even or all odd and $c$ is odd. *Unit net-cubes* are the net-cubes of unit size $\square_{i,j,k} = \square^1_{i,j,k}$.

**Definition 4.1.** A checkerboard IS-net (inscribed spherical net) is a map $f : \mathbb{P} \to \mathbb{R}^3$ satisfying the following conditions:

(i) For any integer $i, j$ the points $\{f_{i,j,k}|k \in \mathbb{Z}\}$ lie on a straight line preserving the order, i.e. the point $f_{i,j,k}$ lies between $f_{i,j,k-1}$ and $f_{i,j,k+1}$. The same holds for points $\{f_{i,j,k}|i \in \mathbb{Z}\}$ and $\{f_{i,j,k}|j \in \mathbb{Z}\}$.

(ii) The unit net-cubes $\square_{i,j,k}$ are circumscribed cubical polytopes, i.e. polyhedra with quadrilateral faces combinatorially equivalent to the three-dimensional cube.

A special case of IS-nets when all cells $\square_{i,j,k}$ with all odd coordinates $i, j, k$ degenerate to points was introduced and investigated by Böhm [10].

An example of an IS-net is shown in Fig. 18.

We denote by $\ell_{j,k}$ the *line of the checkerboard IS-net* that contains the vertexes $\{f_{i,j,k}|\forall i\}$, similarly denote the lines of the other two families by $m_{i,k} \supset \{f_{i,j,k}|\forall j\}$ and $n_{i,j} \supset \{f_{i,j,k}|\forall k\}$. The *planes of the checkerboard IS-net* are denoted by

$$L_i \supset \{f_{i,j,k}|\forall j, k\}, \quad M_j \supset \{f_{i,j,k}|\forall i, k\}, \quad N_k \supset \{f_{i,j,k}|\forall i, j\}.$$

**Theorem 4.1.** *(i) All net-cubes of an IS-net are circumscribed.*

*(ii) The net-cubes* $\square^c_{i,j,k}$ *and* $\square^{c+4s+2}_{i-2s-1,j-2s-1,k-2s-1}$ *are perspective.*

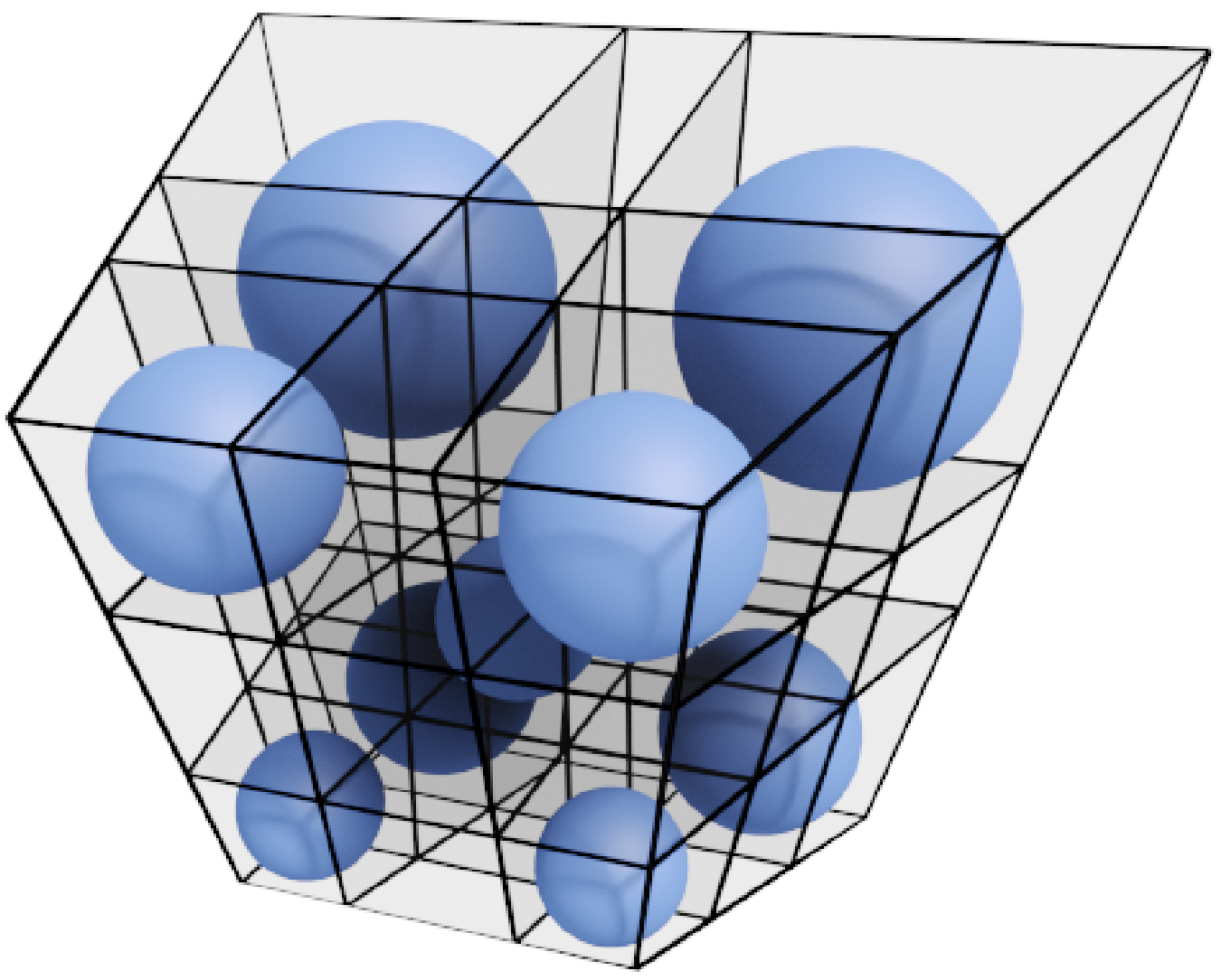

Figure 18: $3 \times 3 \times 3$ block of an IS-net

*(iii)* *For sufficiently large IS-nets (all sides of the cuboid are at least of length 4) all net-cubes are projective images of the standard cube. We will call them* projective cubes*.*

*(iv)* *The lines $\ell_{j,k}$, where $j + k = \mathrm{const}$ lie on a one-sheeted hyperboloid. The same holds for the lines $\ell_{j,k}$ with $j - k = \mathrm{const}$ and for the corresponding lines $m_{i,k}$ and $n_{i,j}$.*

*(v)* *Let $o_{i,j,k}$ be the centers of the spheres inscribed in $\square_{i,j,k}$. The points $o_{i,j,k}$ with all $i, j, k$ even build a grid projectively equivalent to an orthogonal grid (which is built by the intersection points of planes parallel to the coordinate planes). The same claim holds for $o_{i,j,k}$ with all $i, j, k$ odd.*

## 4.2 Construction of checkerboard IS-nets

Our construction of checkerboard IS-nets is based on an incidence theorem on circumscribed cubical polytopes. We start with a rather obvious result.

**Lemma 4.2.** *Given all but one face planes of a circumscribed projective cube, the last face plane is uniquely determined.*

There are three "infinity" points associated with a three-dimensional projective cube. If one deletes two opposite faces of a three-dimensional projective cube the remaining four face plains intersect in a point. The last face plane in the Lemma is uniquely determined by the conditions that it passes through two such "infinity" points and touches the inscribed sphere.

**Lemma 4.3.** *Suppose a cubical polytope* $\square$ *in* $\mathbb{R}^3$ *is split by* 6 *planes in* $27 = 3 \times 3 \times 3$ *combinatorial cubes as shown in Fig. 18. Let us label them naturally by* $\square_{0,0,0}, \square_{1,0,0}, \ldots, \square_{2,2,2}$. *Assume that the central cube* $\square_{1,1,1}$ *is circumscribed and the "frame" cubes* $\square_{0,0,0}$, $\square_{2,0,0}$, $\square_{0,2,0}$, $\square_{0,0,2}$ *are projective cubes and circumscribed as well. Then* $\square$ *is an IS-net* $\square^3_{0,0,0}$, *i.e. the unitary net cubes* $\square_{2,2,0}$, $\square_{2,0,2}$, $\square_{0,2,2}$, $\square_{2,2,2}$ *are circumscribed as well. Moreover the cubes* $\square_{k,l,m}$, $k, l, m \in \{0, 2\}$ *are projective cubes.*

*Proof.* Let us denote by $o_{i,j,k}$ the centers of the spheres inscribed in $\square_{i,j,k}$ (if they exist). Then the central cube $\square_{1,1,1}$ has vertices $f_{1,1,1}$, $f_{2,1,1}$, $f_{1,2,1}$, $f_{1,1,2}$, $f_{2,2,1}$, $f_{1,2,2}$, $f_{2,1,2}$, $f_{2,2,2}$ and its insphere is centered at $o_{1,1,1}$. Consider the projective map $\sigma$ that preserves the point $o_{1,1,1}$ and maps four vertices of the central cube to the centers of the corresponding "frame" spheres: $\sigma(f_{1,1,1}) = o_{0,0,0}$, $\sigma(f_{2,1,1}) = o_{2,0,0}$, $\sigma(f_{1,2,1}) = o_{0,2,0}$, $\sigma(f_{1,1,2}) = o_{0,0,2}$. This projective map preserves four straight lines in general position passing through $o_{1,1,1}$ and therefore preserves all straight lines through $o_{1,1,1}$.

The cells $\square_{2,0,0}$ and $\square_{0,2,0}$ are projective cubes with inscribed spheres. The face planes of the cell $\square_{2,2,0}$ coincide with the corresponding face planes of the cells $\square_{2,0,0}$ and $\square_{0,2,0}$. It it easy to see that this implies that the cell $\square_{2,2,0}$ is circumscribed. Its center is the image of the corresponding vertex $o_{2,2,0} = \sigma(f_{2,2,1})$. Indeed, the plane $L_1$ (see Fig. 19) is mapped by $\sigma$ to the plane $L_{\frac{1}{2}}$ passing through $o_{0,0,0}$, $o_{0,2,0}$ and $o_{0,0,2}$. This plane is the bisector of the planes $L_0$ and $L_1$. The intersection point $L_{\frac{1}{2}} \cap (f_{1,2,2},\, o_{1,1,1}) = \sigma(f_{1,2,2})$ is equidistant from the planes $L_1, M_2, N_2$ as a point of $(f_{1,2,2},\, o_{1,1,1})$ and from $L_0$ as a point of $L_{\frac{1}{2}}$. Thus it is the center of the inscribed sphere. The same argument shows that the cells $\square_{2,0,2}$, $\square_{0,2,2}$ have inscribed spheres centered at $o_{2,0,2} = \sigma(f_{2,1,2})$ and $o_{0,2,2} = \sigma(f_{1,2,2})$ respectively.

Finally, the last corner cell $\square_{2,2,2}$ also has an inscribed sphere. It is centered at the intersection point of three bisector planes $o_{2,2,2} = L_{2\frac{1}{2}} \cap M_{2\frac{1}{2}} \cap N_{2\frac{1}{2}}$. Indeed, this point lies on the straight line $(o_{1,1,1},\, f_{2,2,2})$, and thus is equidistant from all the face planes of $\square_{2,2,2}$.

□

**Corollary 4.4** (9 inspheres incidence theorem)**.** *Suppose a cubical polytope in* $\mathbb{R}^3$ *is split by* 6 *planes in* $27 = 3 \times 3 \times 3$ *combinatorial cubes. Suppose the central and seven of the corner cells are circumscribed. Then the last corner cell is also circumscribed.*

This claim is an immediate corollary of Lemma 4.3. For example, the ‘ìnfinity" points of the projective cube $\square_{2,0,0}$ are homothetic centers of the pairs of spheres inscribed in $\square_{2,0,0}$ and in the cells $\square_{0,0,0}$, $\square_{2,2,0}$ and $\square_{2,0,2}$.

**Theorem 4.5** (global construction of an IS-net from a $2 \times 2 \times 2$ block)**.** *An IS-net is uniquely determined by its* $2 \times 2 \times 2$ *block, which is the union of* $8 = 2 \times 2 \times 2$ *cubic cells with the vertices* $f_{i,j}$, $i, j \in \{0, 1, 2\}$, *and the cells* $\square_{0,0,0}$ *and* $\square_{1,1,1}$ *are circumscribed projective cubes.*

*Proof.* The initial $2 \times 2 \times 2$ block determines uniquely the blocks along the coordinate axes, i.e. all the vertices $f_{i,j,k}$ with the indexes $\{j, k \in \{0, 1, 2\}, \forall i\}$, $\{i, k \in \{0, 1, 2\}, \forall j\}$, $\{i, j \in \{0, 1, 2\}, \forall k\}$. Indeed, there is a unique sphere touching the planes $M_0, M_1, N_0, N_1$ and $L_2$. The corresponding cubic cell should be a circumscribed projective cube. Due to Lemma 4.2 its last face plane $L_3$ is uniquely determined. Proceeding further this way we determine all the planes of the net. Now, consequently applying Lemma 4.3 we prove that all unit net-cubes $\square_{i,j,k}$, $(i, j, k$ are either all even or all odd) are circumscribed. Here we start with the sequence of the cells $\square^3_{0,0,0}$, $\square^3_{1,1,1}$, $\square^3_{2,0,0}$, $\square^3_{3,1,1}$, $\square^3_{4,0,0}$, $\square^3_{5,1,1}$, $\square^3_{6,0,0}$, $\square^3_{7,1,1}$, $\ldots$, and proceed further with the shifted ones $\square^3_{0,2,0}$, $\square^3_{1,3,1}$, $\square^3_{2,2,0}$, $\square^3_{3,3,1}$, $\square^3_{4,2,0}$, $\square^3_{5,3,1}$, $\square^3_{6,2,0}$, $\square^3_{7,3,1}$, $\ldots$, $\square^3_{0,0,2}$, $\square^3_{1,1,3}$, $\square^3_{2,0,2}$, $\square^3_{3,1,3}$, $\square^3_{4,0,2}$, $\square^3_{5,1,3}$, $\square^3_{6,0,2}$, $\square^3_{7,1,3}$, $\ldots$, and so on.

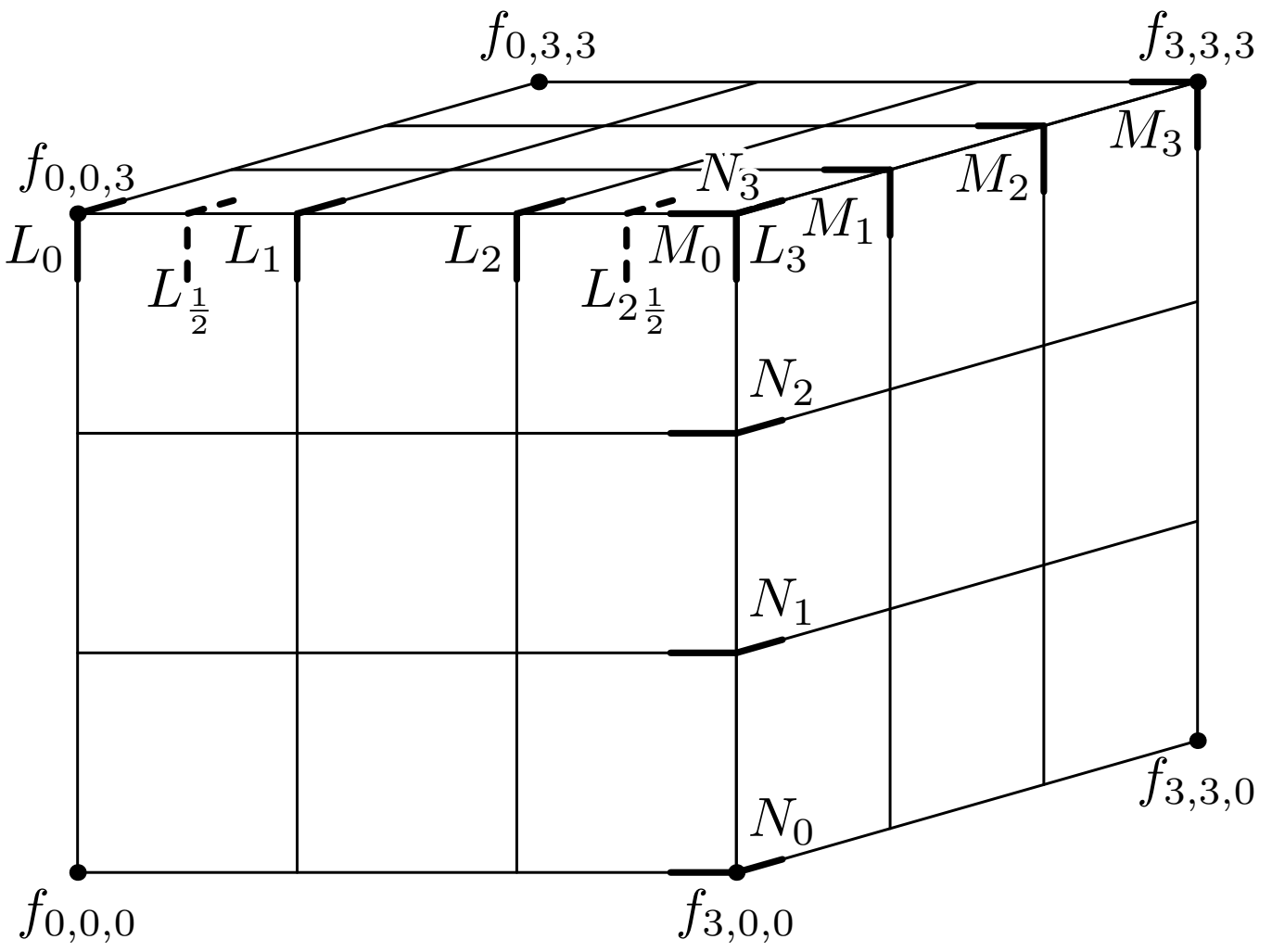

Figure 19: Combinatorics of the IS-net planes. The plane $L_1$ contains the points $f_{1,i,j}$, $\forall i, j$. The planes $L_{\frac{1}{2}}$ and $L_{2\frac{1}{2}}$ are the bisectors of $L_0$ and $L_1$, and of $L_2$ and $L_3$ respectively.

$\square$

## 4.3 Proof of Theorem 4.1

We start with two claims of independent interest. As in previous sections the claims and the proofs can be directly generalized for higher dimensions.

**Lemma 4.6.** *Suppose $\sigma$ is a projective map of $\mathbb{R}^3$ (with infinite plane) to itself which preserves a point $o$ and all lines that pass through $o$. Then the image of any sphere with center at $o$ is a quadric with focus at $o$, i.e. this quadric is an image of rotation of a conic with focus at $o$ around its major axis.*

*Proof.* Consider $\mathbb{R}^3$ embedded in 3-dimensional projective complex space $\mathbb{CP}^3$: the points of $\mathbb{R}^3$ have coordinates $(z_1, \ldots, z_4)$ with $z_4 = 1$ and real $z_1$, $z_2$, and $z_3$. The map $\sigma$ can be naturally extended to a projective transformation of $\mathbb{CP}^3$.

The center of any sphere $\omega \subset \mathbb{CP}^3$ is the pole of the infinite plane $P_\infty$ with respect to $\omega$. The sphere $\omega$ intersects $P_\infty$ in a conic $C_\infty$ given by $\sum_{i=1}^3 z_i^2 = 0$, and any quadric containing $C_\infty$ is a sphere.

Let $o$ be the center of $\omega$. For any $z \in C_\infty$ the line $oz$ touches $\omega$, and all these lines form a cone $C_o$ with the tip $o$ and circumscribed around $\omega$. The cone $C_o$ is fixed under $\sigma$, therefore $\sigma(\omega)$ is also inscribed in it, that means it touches all forming lines of $C_o$.

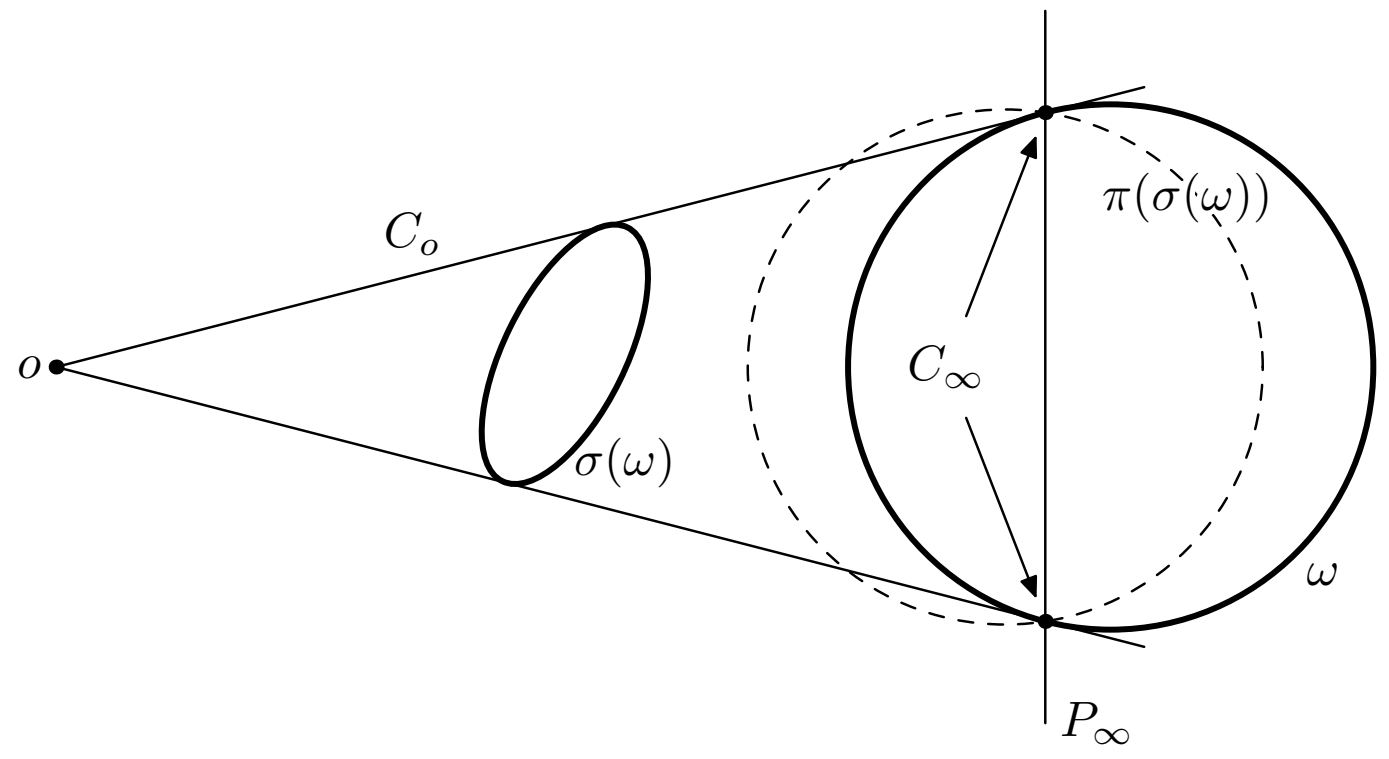

Figure 20: To the proof of Theorem 4.6

In fact, that implies that $o$ is a focus of the real part of $\sigma(\omega)$. Indeed, let $\pi$ be the duality map with respect to the unit sphere $\Omega$ centered at $o$. Note that $\pi(C_o) = C_\infty$, since $\Omega$ touches $C_o$ exactly at $C_\infty$. Since the quadric $\pi(\sigma(\omega))$ contains $C_\infty$ it is a sphere.

Denote by $\omega_\mathbb{R}$ the points of $\omega$ with real coordinates. We know that $\pi(\sigma(\omega_\mathbb{R}))$ is a sphere.

That implies that $\sigma(\omega_\mathbb{R}) = \pi(\pi(\sigma(\omega_\mathbb{R})))$ is an image of rotation of a conic with focus at $o$ around its major axis. Indeed, it is known that the polar image of a circle with respect to another circle is a conic with a focus at the center of the latter circle (see [4]). Therefore, for any plane $L$ passing through $o$ and the center of the sphere $\omega_\mathbb{R}$ we have:

$$\sigma(\omega_\mathbb{R}) \cap L = \pi(\pi(\sigma(\omega_\mathbb{R}))) \cap L = \pi_L(\pi_L(\sigma(\omega_\mathbb{R} \cap L))) = \gamma_L.$$

Here $\pi_L$ is the two-dimensional polar transformation with respect to the unit circle lying in $L$ and centered at $o$, and $\gamma_L$ is a conic in $L$ with a focus at $o$. Finally, the rotational symmetry in the choice of $L$ implies the claim. $\square$

**Theorem 4.7.** *Suppose a cubical polytope ⧉ in $\mathbb{R}^3$ is split by 6 hyperplanes in $27 = 3 \times 3 \times 3$ cubical polytopes. Suppose the central and $8 = 2 \times 2 \times 2$ corner cells are circumscribed. Then the cube ⧉ is also circumscribed.*

*Proof.* Let $J$, $|J| = 8$, be a set of indices of corner cells ⧉$_j$, $j \in J$. Denote by $o_j$ the centers of their inscribed spheres, and by $a_j$ and $b_j$ the corresponding vertices of ⧉ and of the central cell respectively. Let $o$ be the center of the inscribed sphere of the central cell.

As was mentioned in the proof of Lemma 4.3, there is a projective map $\sigma$ which preserves the point $o$ and maps points $b_j$ to points $o_j$. Denote the image of the central cell under $\sigma$ by ⧉$'$.

From Lemma 4.6 it follows that $\sigma$ maps the inscribed sphere of the central cubical polytope to an ellipsoid $\alpha$ with focus at $o$. This ellipsoid touches faces of ⧉$'$. Let us show that the second focus $o''$ of this ellipsoid is the center of the sphere inscribed in $\square$.

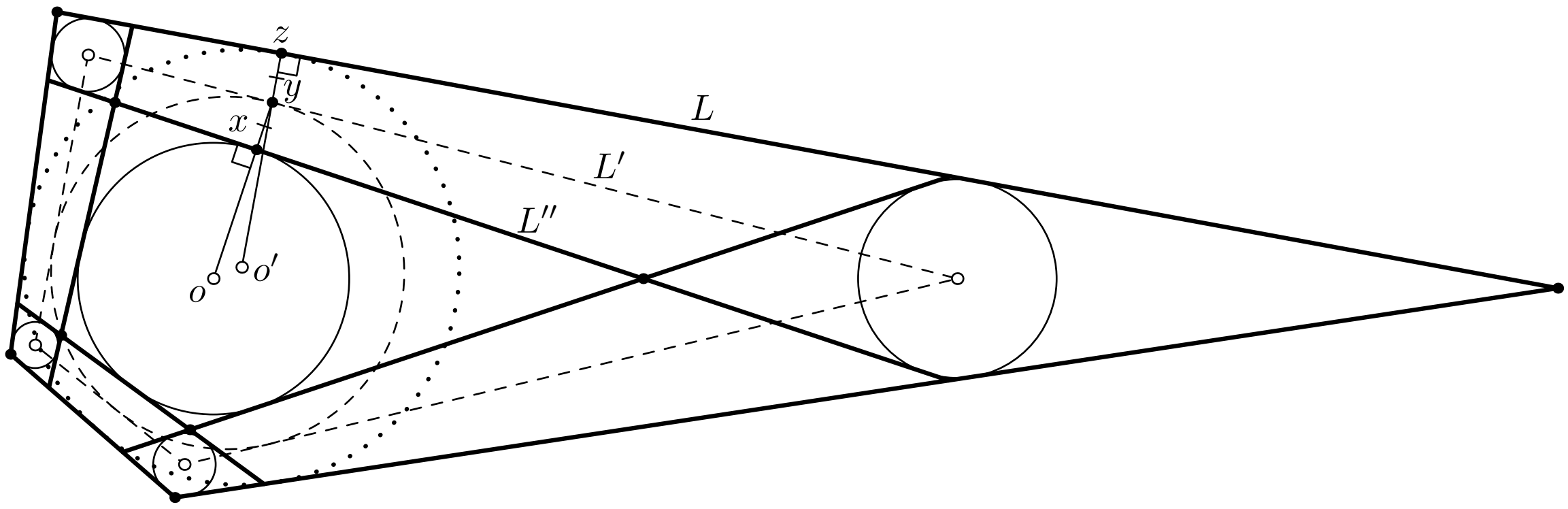

Figure 21: To the proof of Theorem 4.7

Take a face plane $L$ of $\square$, let $L'$ and $L''$ be the corresponding face planes of $\square'$ and of the central cell respectively. Let the inscribed sphere of the central cell touches $L''$ at $x$, then the ellipsoid $\alpha$ touches $L'$ at a point $y$ lying on the line $(ox)$. Let $z$ be the point symmetric to $x$ with respect to the hyperplane $L'$ (the two-dimensional analogue of this construction is shown on Fig. 21). The optical property of ellipsoids implies that $z$, $y$ and $o'$ are collinear. Note that the corresponding line is perpendicular to $L$ since the line $(xy)$ is perpendicular to $L''$. Thus,

$$|zo'| = |zy| + |yo'| = |xy| + |yo'| = l - |ox|,$$

where $l$ is the big axis of the ellipsoid $\alpha$. Since $l - |ox|$ is independent of the choice of the $L$, we have that $o'$ lies at equal distances from all face planes of $\square$. □

*Proof of Theorem 4.1.* (i) follows from Theorem 4.7 which we apply inductively similarly to the proof of (i) of Theorem 3.1.

(ii) Let $s$ be positive (negative case proved analogously). Sides of the small cube $\square^c_{i,j,k}$ divide the big cube $\square^{c+4s+2}_{i-2s-1,j-2s-1,k-2s-1}$ into 27 cells, note that 8 corner cells are net-cubes. Denote them by $\square_l$ where $l = 1, \dots, 8$. Denote by $a_l$ and by $b_l$ the corresponded vertices of $\square^c_{i,j,k}$ and $\square^{c+4s+2}_{i-2s-1,j-2s-1,k-2s-1}$ respectively. Let $\omega_l$ be the inscribed spheres of $\square_l$ and $\Omega$ and $\Omega'$ be inscribed spheres of $\square^c_{i,j,k}$ and $\square^{c+4s+2}_{i-2s-1,j-2s-1,k-2s-1}$.

Note that $a_l$ is the negative homothety center of $\Omega'$ and $\omega_l$. Analogously $b_l$ is the positive homothety center of $\Omega$ and $\omega_l$. From the Monge theorem it follows that the line $(a_l b_l)$ passes through the negative homothety center of $\Omega$ and $\Omega'$. Therefore all lines $(a_l b_l)$ pass through one point.

(iii) Consider the net-cube $\square^c_{i,j,k}$. The planes $L_i$, $L_{i+c}$, $M_j$, and $M_{j+c}$ pass through the center of the positive homothety of spheres inscribed in $\square^c_{i,j,k}$ and $\square^c_{i,j,k+2}$. This is one of three ‘ìnfinity” points of a projective cube. The existence of other two ‘ìnfinity” points is proven in the same way.

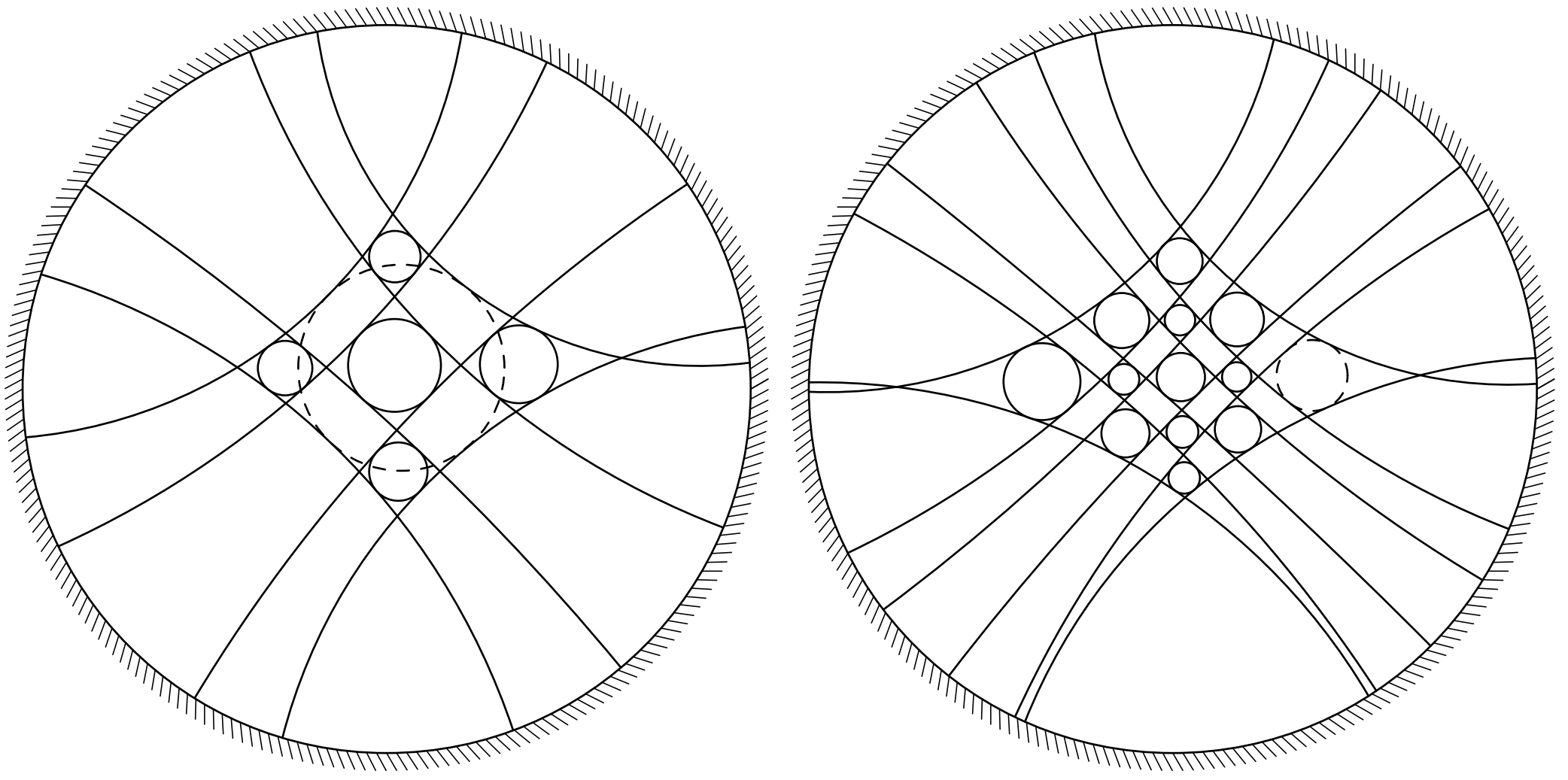

Figure 22: Hyperbolic six incircles lemma and checkerboard incircles incidence theorem

(iv) Let $\ell_{i_1,j_1}$ and $\ell_{i_2,j_2}$ be two lines where $i_1$ and $i_2$ are of different parities and $i_1 + j_1 = i_2 + j_2$. Then there is a net-cube with edges on $\ell_{i_1,j_1}$ and $\ell_{i_2,j_2}$. Since net-cubes are projective cubes these lines intersect.

So we can separate the family of lines $\ell_{i,j}$ into two subfamilies with odd and even $i$. Any two lines from different families intersect (or parallel). Therefore they lie on a hyperboloid of one sheet and these lines are asymptotic lines of two families on the hyperboloid.

(v) As in the proof of Lemma 4.7 we show that the points $o_{i+1\pm1.j+1\pm1,k+1\pm1}$ are vertices of a cube projectively equivalent to $\square_{i+1,j+1,k+1}$. The latter is a projective cube due to (iii). Thus all elementary cells of the grid $o_{i,j,k},\ \ i,j,k \in 2\mathbb{Z}$ are projective cubes. It is easy to see that the infinite points of all elementary cells coincide, which implies the claim.

□

# 5 IC- and IS-nets in hyperbolic space

For almost all results of the previous sections there exist natural analogues in hyperbolic and spherical spaces.

In this section we discuss the case of hyperbolic space. For simplicity we assume that the combinatorics of the corresponding geodesics coincides with those in Euclidean space. In particular, we assume that the lines from different coordinate families intersect. In this case the embedded pieces $f : \mathbb{P} \to H$ of patterns in the hyperbolic space look similar to their Euclidean analogues, and the claims and the proofs almost coincide.

*On Section 2.*

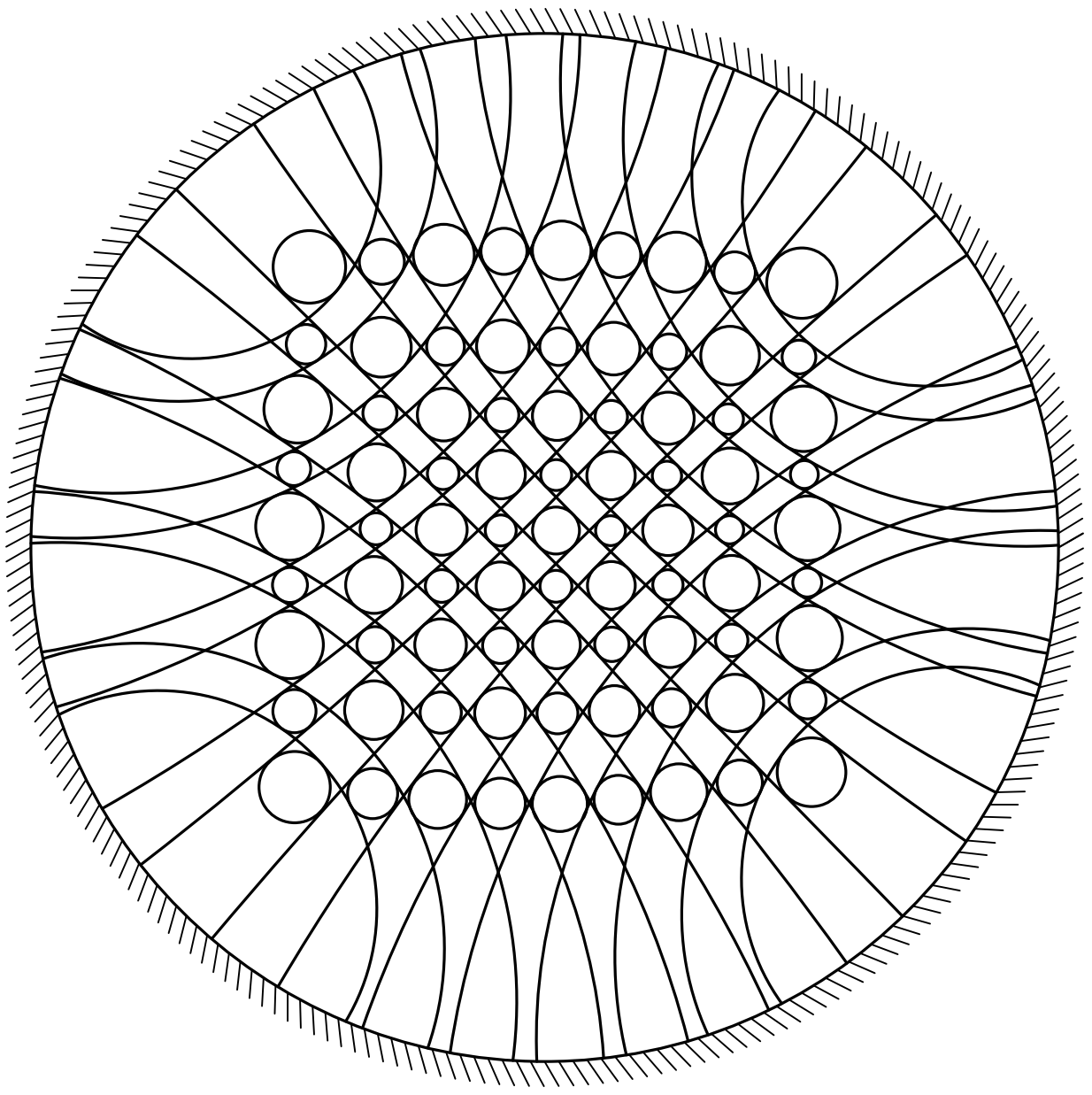

Figure 23: Hyperbolic checkerboard IC-net

For the proof of hyperbolic version of Theorem 2.5 we use the Klein model. It is well-known that conics (as a set such that sum or difference of distances to two fixed points is constant) in the Klein model are Euclidean conics, and confocal conics build a dual pencil containing the circle of absolute as an element [16]. Thus, we can use the Euclidean version of Lemma 2.4 for the proof, as we did for Theorem 2.5. In that proof we also use the equal angle lemma, which is valid in the hyperbolic plane as well.

Corollaries and Lemmas 2.6-2.8 have similar proofs for hyperbolic case as well as Theorem 2.9, which follows from them.

Since in the hyperbolic space there is no homothety transformation, IC-nets build a five-dimensional family. The construction of Corollary 2.10 is valid in the hyperbolic plane as well.

Items (i)-(iv), (viii), and (ix) of Theorem 2.1 hold true in the hyperbolic case as well. The proof is the same, except of (ix). Here we observe that the picture in the Klein model coincides with the Euclidean one: the bisectors are tangent to the conic $\alpha'$ and their points of intersection lie on the conic dual to $\alpha$ with respect to $\alpha'$.

*On Section 3.*

All claims of Theorem 3.1 except (vi) hold true also in the hyperbolic case. In the proof we used the following results, which are all valid in the hyperbolic case as well:
Lemma 3.2;
the Pascal theorem (conics in Klein model are Euclidian conics and geodesics are straight lines);
the Monge theorem (it can be proven through a special construction in three-dimensional space, which works also in hyperbolic plane (see [2]).

To prove the hyperbolic version of Theorem 3.3 (see Fig. 22, right) we repeat the construction of Section 3.2. Consider the conformal ball model of the three dimensional hyperbolic space $H^3$, and take its equator plane as a conformal model of the two dimensional hyperbolic space $H^2$. Draw the circles and geodesics pattern in $H^2$ as shown in Fig. 13. Construct the corresponding cones intersecting $H^2$ at constant angle. As in the Euclidean case of Theorem 3.3 the corresponding apexes are collinear. Passing to the Klein model we obtain a projective picture, and the projective incidence Theorem 3.4 completes the proof.

Checkerboard confocal IC-nets in the hyperbolic space are defined exactly as in subsection 3.3 and a hyperbolic version of Theorem 3.7 holds.

*On Section 4.*

In hyperbolic space Theorem 4.1 holds with a modified version of items (iii) and (v). A "projective cube" there means a polytope which in the Klein models is projective image of a cube. The proof is the same as for Euclidean case, but hyperbolic versions of lemmas used in the proof require some comments.

The proof of Lemma 4.3 consists of two parts: an observation that the centers of certain inscribed spheres lie on planes and an incidence theorem of lines and planes, which we prove by considering appropriate projective transformation. In the hyperbolic case we have the same condition on centers of inscribed spheres and the required incidence theorem in the Klein model is equivalent to the Euclidean one.

For the proof of the hyperbolic version of Theorem 4.7 we use Lemma 4.6 with the following observation: if a $o$ focus of a quadric coincides with center of the Klein model then it is also an Euclidean focus of it (and vise versa). We move $o$ to the center of the Klein model, by Lemma 4.6 there is an (Euclidean) ellipsoid touching $\square'$ with focus at $o$, and thus its hyperbolic focus is also at $o$. Now we choose the second hyperbolic focus $o'$ and repeat the Euclidean arguments.

# 6 Appendix. A direct proof of the Graves–Chasles theorem

Here we give a direct computational proof of the Graves–Chasles Theorem 2.5.

*Proof.* Let us show that (ii) implies (i). We assume that $\alpha$ is an ellipse with foci $f_1$ and $f_2$. The case of hyperbola can be proven in the same way. Let $\alpha_1$ be the ellipse passing through the points $a$ and $c$, and $b_1 := (af_1) \cap (bf_2)$, $d_1 := (af_2) \cap (bf_1)$. Let $f_1'$ be a reflection of $f_1$ in line $(ab)$. Note that $|f_1'f_2|$ equals the length $l$ of the the big axis of the ellipse $\alpha$. We denote by $l_1$ the length of the big axis of $\alpha_1$.

Since $a$ and $c$ lie on the ellipse $\alpha_1$, we have $|f_1a|+|f_2a| = |f_1c|+|f_2c|$. Therefore the quadrilateral $(ab_1cd_1)$ is circumscribed. Denote the center of its circle by $o$. From Poncelet's isogonal lemma (see, for example, Theorem 1.4 in [4]) it follows that $\angle b_1ab = \angle dad_1$ and $\angle b_1cb = \angle dcd_1$. So, it is sufficient to show that the distances from $o$ to $(ad)$ and $(cd)$ are equal. We observe that

$$\frac{d(o,ad)}{d(o,ad_1)} = \frac{\sin\angle iad}{\sin\angle iad_1} = \frac{\cos\frac{1}{2}(\pi-\angle bad)}{\cos\frac{1}{2}(\pi-\angle b_1ad_1)} = \frac{\cos\frac{1}{2}\angle f_1'af_2}{\cos\frac{1}{2}\angle f_1af_2}.$$

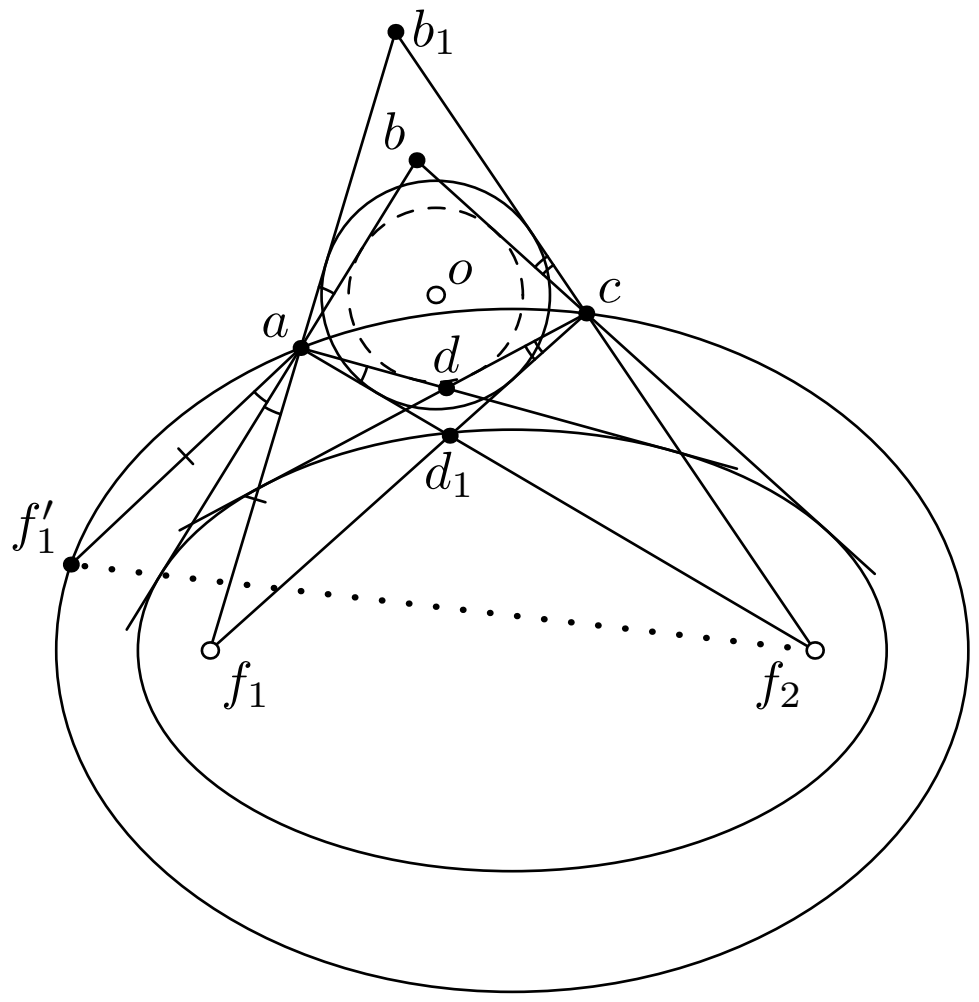

Figure 24: To the proof of the Graves–Chasles theorem

Further $2(\cos\frac{1}{2}\angle f_1'af_2)^2 = \cos\angle f_1'af_2 + 1$ holds, and from the cosine formula we obtain $\cos\angle f_1'af_2 = \frac{|f_1a|^2+|f_2a|^2-l^2}{2|f_1a|\cdot|f_2a|}$. And we get $\cos\angle f_1'af_2/2 = \sqrt{\frac{l_1^2-l^2}{2|f_1a|\cdot|f_2a|}}$. An analogous computation gives $2(\cos\frac{1}{2}\angle f_1af_2)^2 = \cos\angle f_1af_2 + 1$, $\cos\angle f_1af_2 = \frac{|f_1a|^2+|f_2a|^2-|f_1f_2|^2}{2|f_1a|\cdot|f_2a|}$, $\cos\frac{1}{2}\angle f_1af_2 = \sqrt{\frac{l_1^2-|f_1f_2|^2}{2|f_1a|\cdot|f_2a|}}$. Therefore

$$\frac{d(o,ad)}{d(o,ad_1)} = \frac{\cos\frac{1}{2}\angle f_1'af_2}{\cos\frac{1}{2}\angle f_1af_2} = \sqrt{\frac{l_1^2-l^2}{l_1^2-|f_1f_2|^2}}.$$

We see that the ratio $\frac{d(o,ad)}{d(o,ad_1)}$ is independent of the point $a$. Hence

$$\frac{d(o,ad)}{d(o,ad_1)} = \frac{d(o,bd)}{d(o,bd_1)},$$

and finally $d(o,ad) = d(o,cd)$, since $d(o,ad_1) = d(o,cd_1)$.

Let us show now that (i) implies (ii).

If $c$ does not lie on the ellipse with foci $f_1$ and $f_2$ passing through $a$ we can choose another point $c'$ on $(bc)$ such that it is and define $d'$ as the point of intersection of $(ad)$ with the tangent line from $c'$ to $\alpha$. The quadrilateral $abc'd'$ is circumscribed. But on the other hand, the incircles of $(abc'd')$ and $(abcd)$ coincide, and $(cd)$ and $(c'd')$ are the common interior tangent lines of $\alpha$ and this incircle. They should coincide, thus $c_1 = c$. The equivalence (i) $\Leftrightarrow$ (iii) can be shown in the same way. $\square$